\documentclass[12pt, twoside]{article}
\usepackage[pagewise]{lineno}
\usepackage{amsmath,amssymb}
\usepackage{float}
\usepackage{url}
\usepackage{multirow}
\usepackage{latexsym}
\usepackage{bm}
\usepackage{graphicx}
\graphicspath{{fig_eps/}}
\usepackage{epsfig}
\usepackage{cite}
\usepackage{epstopdf}
\usepackage{caption}
\usepackage{subfigure}
\usepackage[justification=centering]{caption}
\usepackage{textcomp,booktabs}
\usepackage[usenames,dvipsnames]{color}
\usepackage{colortbl}

\definecolor{mygray}{gray}{0.9}
\definecolor{mypink}{rgb}{0.99,0.91,0.95}
\definecolor{mycyan}{cmyk}{0.3,0,0,0}
\usepackage{booktabs}
\usepackage[top=2.5cm,bottom=2.5cm,left=2cm,right=2cm]{geometry}
\usepackage[colorlinks,linkcolor=red,anchorcolor=green,citecolor=blue,CJKbookmarks=True]{hyperref}
\newtheorem{theorem}{Theorem}[section]
\newtheorem{lemma}{Lemma}[section]
\newtheorem{remark}{Remark}[section]
\newtheorem{assumption}{Assumption}[section]
\newtheorem{definition}{Definition}[section]
\newtheorem{proposition}{Proposition}[section]
\allowdisplaybreaks
\def\Proof{\noindent{\bf Proof.}~}
\def\qed{\hfill$\square$\smallskip}

\begin{document}

\title{\textbf{Traveling Waves of Modified Leslie-Gower Predator-prey Systems}}
\author{Hongliang Li\thanks{Corresponding author.}, Min Zhao, Rong Yuan\\
Laboratory of Mathematics and Complex Systems (Ministry of Education)\\
School of Mathematical Sciences, Beijing Normal university\\
Beijing 100875, People's Republic of China\\
E-mail: lihl@mail.bnu.edu.cn,
minzhao@mail.bnu.edu.cn,\\ ryuan@bnu.edu.cn
}

\date{}

\maketitle

\begin{abstract}
The spreading phenomena in modified Leslie-Gower reaction-diffusion predator-prey systems are the topic of this paper. We mainly study the existence of two different types of traveling waves. Be specific, with the aid of the upper and lower solutions method, we establish the existence of traveling wave connecting the prey-present state and the coexistence state or the prey-present state and the prey-free state by constructing different and appropriate Lyapunov functions. Moreover, for traveling wave connecting the prey-present state and the prey-free state, we gain more monotonicity information on wave profile based on the asymptotic behavior at negative infinite. Finally, our results are applied to modified Leslie-Gower system with Holling II type or Lotka-Volterra type,
and then a novel Lyapunov function is constructed for the latter, which further enhances our results. Meanwhile, some numerical simulations are carried to support our results.


\end{abstract}

\section{Introdcution}
\noindent

Due to the complexity of ecosystems, long-term coexisting species frequently exhibit diversified interspecific interactions, such as cooperation, competition, and predation. One of the most famous systems is Lotka-Volterra predator-prey system derived by Lotka \cite{a25} and Volterra \cite{v26} in the 1920s, and there has been a substantial amount of research devoted to predator-prey systems, seeing \cite{ccr09,m89,h66,rm63,ag89,ah00}.
However, the majority of biological species live in spatially heterogeneous natural habitats, and it is fair to anticipate that the variability of the habitats will have an impact on population density. Consequently, the response terms should be expanded to include the diffusion.

What we are curious about is the change in ecological processes when alien species are introduced into new habitats. One way to study this topic is to look at the so-called traveling waves, seeing \cite{D83,D84,hjh03,h12,hsu12,dlm21}, another way is to characterize spreading speeds of the predator by solving the Cauchy problem, seeing \cite{a75,a78,d19}. It should be emphasized that the functional response, the rate of prey consumption by an average predator, and the predator growth are identical in most two species predator-prey systems. For example,
\begin{equation*}
\left\{\begin{split}
&\displaystyle\frac{\partial u(x, t)}{\partial t}=\Delta u(x, t) +g(u(x, t))-f(u(x, t))v(x, t),\\[0.2cm]
&\displaystyle\frac{\partial v(x, t)}{\partial t}=d\Delta v(x, t) +\left(\beta f(u(x, t))-d\right)v(x, t).
\end{split}
\right.
\end{equation*}
While some scholars have also studied predator-prey systems whose predator growth function is different from the functional response. For example, Leslie in \cite{l48,l60} proposed a predator-prey system in which the predator's environmental carrying capacity is proportional to the abundance of the prey
\begin{equation*}
\left\{\begin{split}
&\displaystyle\frac{\partial u(x, t)}{\partial t}=\Delta u(x, t)+u(x, t)\left(1-u(x, t)-\frac{av(x, t)}{u(x, t)+e_1}\right), \\[0.2cm]
&\displaystyle\frac{\partial v(x, t)}{\partial t}=d \Delta v(x, t)+ sv(x, t)\left(1-\frac{v(x, t)}{u(x, t)}\right),
\end{split}
\right.
\end{equation*}
where the term $v(x, t)/u(x, t)$ is called the Leslie-Gower term. Such predators are called \emph{specialist predators} \cite{t13} since the population density of a favored prey is small, the predator should concentrate on other prey species, and vice versa.
Recently, the existence of traveling waves to the above system has attracted increasing interest, including Lotka-Volterra type, Holling-III and Beddington-DeAngelis functional responses and so on, seeing \cite{cgy17,adp17,wf21,D12,zw20}.

Based on the above system, the authors in \cite{a02,a03} proposed a modified Leslie-Gower predator-prey system by adding the positive constant $\mu$ to the Leslie-Gower term. Such predators can be called \emph{generalist predators}\cite{g19} since they explore other alternative food sources in the absence of their favorite prey. More precisely,
\begin{equation}\label{eq:mlg}
\left\{\begin{split}
&\displaystyle\frac{\partial u(x, t)}{\partial t}=\Delta u(x, t)+u(x, t)\left(1-u(x, t)-\frac{av(x, t)}{u(x, t)+e_1}\right), \\[0.2cm]
&\displaystyle\frac{\partial v(x, t)}{\partial t}=d \Delta v(x, t)+sv(x, t)\left(1-\frac{v(x, t)}{u(x, t)+\mu}\right),
\end{split}
\right.
\end{equation}
where parameters $a$,
$s$, $e_1$ and $\mu$ are positive, $u(x, t)$ and $v(x, t)$ represent population densities of prey and predator species respectively, and $d$ denotes the ratio of the diffusion of the predator to that of the prey. For more detailed descriptions, one can refer to \cite{a02,a03}.

When the prey is stationary, Tian et al. in \cite{t16,t18} studied traveling waves connecting two appropriate equilibria by the shoot argument and the upper and lower solutions method respectively. However, when the prey diffuses slower than the predator, Hsu and Lin in \cite{hl19} obtained the existence of traveling waves connecting $(0,0)$ and the coexistence state by the upper and lower solutions method with the aid of contracting rectangles. For system \eqref{eq:mlg}, to the best of the author's knowledge, there are no results on traveling waves connecting the prey-present state and the coexistence state or the prey-present state and the prey-free state, even for a modified Leslie-Gower predator-prey system with Lotka-Volterra type \cite{f80, g12}. Based on these, we consider a generalized Leslie-Gower system
\begin{equation}\label{eq:h11}
\left\{\begin{split}
&\displaystyle\frac{\partial u(x, t)}{\partial t}=\Delta u(x, t)+f(u(x, t))\left(p(u(x, t))-v(x, t)\right), \\[0.2cm]
&\displaystyle\frac{\partial v(x, t)}{\partial t}=d \Delta v(x, t)+sv(x, t)\left(1-\frac{v(x, t)}{u(x,t)h(u(x,t))+\mu}\right).
\end{split}
\right.
\end{equation}
where parameters $d,s,\mu$ are positive constants.
In this paper, we assume that $f(u)$, $p(u)$ and $h(u)$ are $C^1$ functions satisfying the following assumption:
\begin{assumption}\label{as}
	\text{ }\\[0.2cm]
    (H1) $f(0)=0$ and $f'(u)>0$ for $u\in[0,+\infty)$.\\[0.2cm]
    (H2) $p(1)=0$ and $p'(u)<0$ for $u\in(0,+\infty)$.\\[0.2cm]
    (H3) $h(0)>0$ and $h'(u)\geq0$ for $u\in[0,+\infty)$.
\end{assumption}
Some classical examples satisfying Assumption \ref{as} as follows
\begin{description}
  \item (1) Lotka-Volterra type system: for given $a>0$,
  \begin{equation*}
    f(u)=au, \quad p(u)=\displaystyle\frac{(1-u)}{a}\quad \text{ and }\quad h(u)=1.
  \end{equation*}
  \item (2) Holling II type system: for given $a>0$ and $e_1\geq1$,
  \begin{equation*}
    f(u)=\displaystyle\frac{au}{u+e_1}, \quad p(u)=\displaystyle\frac{(1-u)(u+e_1)}{a}\quad \text{ and }\quad h(u)=1.
  \end{equation*}
  \item (3) Ivlev-type type system: for given $a>0$ and $2>m>0$,
  \begin{equation*}
    f(u)=a(1-e^{-mu}), \quad p(u)=\displaystyle\frac{u(1-u)}{a(1-e^{-mu})} \quad \text{ and }\quad h(u)=1.
  \end{equation*}
\end{description}

In view of Assumption \ref{as}, as for the corresponding kinetic system, there are always three boundary equilibria $(0,0)$, $(1,0)$ and $(0,\mu)$. Clearly, $(1,0)$ is unstable. Moreover, one can readily verify that $(0,\mu)$ is stable if $\mu>p(0)$, whereas $(0,\mu)$ is unstable if $\mu<p(0)$. On the other hand, we also conclude that function $Q(u)=p(u)-uh(u)-\mu$ satisfies $Q(1)<0$, $Q(0)>0$ and $Q'(u)<0$ for $u\in(0,1)$ if $p(0)>\mu$, which implies that system \eqref{eq:h11} has a unique positive equilibrium $(u^*,v^*):=(u^*,u^*h(u^*)+\mu)$ with $u^*\in(0,1)$. Consequently, we may expect traveling waves connecting $(1,0)$ and $E$ of system \eqref{eq:h11}, where
\begin{align*}
  E=\left\{
  \begin{array}{ll}
    (0,\mu), &if\ \mu>p(0), \\[0.2cm]
    (u^*,v^*),\ &if\ \mu<p(0).
  \end{array}
  \right.
\end{align*}

We are interested in the spreading phenomena of system \eqref{eq:h11}. First of all, let us focus on spreading speeds of $v(x, t)$. Constant $s^*$ is called \emph{spreading speeds} of species $w(x,t)$ if
\begin{align*}
  &\lim\limits_{t\rightarrow\infty}\left\{\sup\limits_{|x|>(s^*+\epsilon)t}w(x,t)\right\}=0, \text{ for }  \epsilon>0,\\[0.2cm]
  &\liminf\limits_{t\rightarrow\infty}\left\{\inf\limits_{|x|<(s^*-\epsilon)t}w(x,t)\right\}>0, \text{ for } \epsilon\in(0,s^*).
\end{align*}
We assume that $u(x,0)=1$ and $v(x,0)\in[0,h(1)+\mu]$ is continuous function with nonempty compact support, by the theory of reaction-diffusion system, then $[0,1]\times[0,h(1)+\mu]$ is a invariant region of system \eqref{eq:h11}. Hence, we obtain
$(u(x,t),v(x,t))\in[0,1]\times[0,h(1)+\mu]$ for $x\in\mathbb{R}$ and $t>0$.
Meanwhile, one can directly observe that
\begin{equation*}
\displaystyle\frac{\partial v(x, t)}{\partial t}\leq d \Delta v(x, t)+sv(x, t)\left(1-\frac{v(x, t)}{h(1)+\mu}\right),
\end{equation*}
and
\begin{equation*}
\displaystyle\frac{\partial v(x, t)}{\partial t}\geq d \Delta v(x, t)+sv(x, t)\left(1-\frac{v(x, t)}{\mu}\right),
\end{equation*}
thus, the spreading theory and the comparison principle of \cite{a75,a78} give
\begin{align*}
  &\lim\limits_{t\rightarrow\infty}\left\{\sup\limits_{|x|>(2\sqrt{ds}+\epsilon)t}v(x,t)\right\}=0, \text{ for } \epsilon>0,\\[0.2cm]
  &\liminf\limits_{t\rightarrow\infty}\left\{\inf\limits_{|x|<(2\sqrt{ds}-\epsilon)t}v(x,t)\right\}>\mu, \text{ for } \epsilon\in(0,2\sqrt{ds}),
\end{align*}
which implies that spreading speeds of $v(x,t)$ is $c^*=2\sqrt{ds}$. Therefore, in this paper, we mainly establish the existence of traveling waves connecting $(1,0)$ and $E$ of system \eqref{eq:h11}. A positive solution is called a traveling wave of system \eqref{eq:h11} if it has the form
\begin{equation*}
  (u(x,t),v(x,t))=(\tilde u(z),\tilde v(z)), \ z=x+ct,
\end{equation*}
where $c$ is the wave speed. With the tilde removed, $(u(z),v(z))$ should satisfy
\begin{equation}\label{eq:tr}
\left\{\begin{split}
&u''(z)-cu'(z)+f(u(z))\left(p(u(z))-v(z)\right)=0,\\[0.2cm]
&dv''(z)-cv'(z)+sv(z)\left(1-\frac{v(z)}{u(z)h(u(z))+\mu}\right)=0.
\end{split}
\right.
\end{equation}
Further, if $(u(z),v(z))$ also satisfies the boundary condition
\begin{equation}\label{eq:b}
  (u(-\infty),v(-\infty))=(1,0),\quad (u(+\infty),v(+\infty))=E,
\end{equation}
then it is referred to as \emph{strong-traveling wave} of system \eqref{eq:h11}. In addition,
one that just meets the left-hand tail limit is referred to as \emph{semi-traveling wave} of system \eqref{eq:h11}.

In the sequence, strong or semi-traveling waves as mentioned are the above waves. Moreover, for the sake of convenience, we assume that
\begin{equation*}
  q(u)=uh(u)+\mu,
\end{equation*}
clearly, $q(0)=\mu$, $q(1)=h(1)+\mu$ and $q(u)$ is monotone increasing function in $[0,+\infty)$.

We are now in the position to state the existence results. The first theorem states necessary and sufficient condition on the existence of semi-traveling wave.
\begin{theorem}\label{th:ma}
Assume that (H1)-(H3) hold. System \eqref{eq:h11} has a semi-traveling wave if and only if $c\geq c^*:=2\sqrt{ds}$.
\end{theorem}

The following theorem establishes the existence of traveling wave connecting $(1,0)$ and $(u^*,v^*)$ of a class of modified Leslie-Gower predator-prey systems, that is $h(u)=1$.
\begin{theorem}\label{th:ma2}
Assume that (H1)-(H2) hold and $h(u)=1$. For $c\geq c^*$, if  $p(0)>\mu$ and
\begin{description}
  \item (P) $(u^*+\mu)f(u)-(u-u^*)(u+\mu)f'(u)>0$ for $u\in(0,1)$,
\end{description}
then system \eqref{eq:h11} has a traveling wave connecting $(1,0)$ and $(u^*,v^*)$.
\end{theorem}

Finally, we state the existence result on traveling wave connecting $(1,0)$ and $(0,\mu)$. Before we can do that, we define the sets
\begin{align*}
  \mathcal{A}=&\left\{v\in C^2\left(\mathbb{R}\right)\mid 0<v(z)<\mu \ \text{and}\ v'(z)>0 \ \text{over}\ \mathbb{R}\right\}, \\[0.2cm]
  \mathcal{B}=&\left\{
           v\in C^2\left(\mathbb{R}\right) \bigg|
           \begin{array}{ll}
              \exists \ z_v\  \text{such that}\ v(z_v)=\mu, v'(z)>0\ \text{in}\ (-\infty, z_v) \\
              \text{and}\ v(z)>\mu\ \text{in}\ (z_v,+\infty)
           \end{array}
           \right\}.
\end{align*}
\begin{theorem}\label{th:ma3}
Assume that (H1)-(H3) hold.  For $c\geq c^*$, if
\begin{equation*}
  f(u)=ug(u) \quad \text{and} \quad p(u)=(1-u)/g(u),
\end{equation*}
where $g(u)$ satisfies
  \begin{equation*}
    \min_{u\in[0,1]}g(u)>\frac{h(1)(h(1)+\mu)}{\mu^2h(0)},
  \end{equation*}
  then the following hold.
  \begin{description}
    \item (a) System \eqref{eq:h11} has a traveling wave  $(u(z),v(z))$ connecting $(1,0)$ and $(0,\mu)$.
    \item (b) $v(z)\in \mathcal{A}\cup \mathcal{B}$ and $u'(z)<0$ over $\mathbb{R}$.
  \end{description}
\end{theorem}

This paper is organized as follows. In Section 2, we introduce some preliminary knowledges.
In Section 3, based on the upper and lower solutions method and the unstable manifold theorem, we obtain necessary and sufficient condition on the existence of semi-traveling wave. In  section 4, we establish that the existence of strong-traveling waves, and for traveling wave connecting $(1,0)$ and $(0,\mu)$, we further gain more monotonicity information on wave profile based on the asymptotic behavior at $z=-\infty$. In the last section, we apply Theorems \ref{th:ma2} and \ref{th:ma3} to modified Leslie-Gower system with Holling II type or Lotka-Volterra type, then a novel Lyapunov function is constructed for the latter, which further enhances our results. Moreover, the numerical simulations of traveling waves are carried to support our results.

\section{Preliminary}
\noindent

In this section, we will introduce the upper and lower solutions method, which actually
provides a constructive approach for producing traveling waves.

Let us begin with the definition of upper and lower solutions of system \eqref{eq:tr}. Define sets
\begin{equation*}
  \begin{array}{l}
X=\left\{\Phi(z)=(u(z), v(z)) \mid \Phi(z) \in C\left(\mathbb{R}, \mathbb{R}^{2}\right)\right\}, \\[0.3cm]
X_{0}=\left\{\Phi(z)\in X \mid 0\leq u(z)\leq 1,0\leq v(z)\leq q(1), z \in \mathbb{R}\right\}.
\end{array}
\end{equation*}
\begin{definition}\label{Def}
Function pairs $\left(\overline{u}(z),\overline{v}(z)\right)$ and $\left(\underline{u}(z),\underline{v}(z)\right)$ in $X_0$ are called a pair of upper and lower solutions of  system \eqref{eq:tr} if they satisfy
\begin{description}
  \item(a)\  $\underline{u}(z)\leq\overline{u}(z)$, $\underline{v}(z)\leq\overline{v}(z)$ for $z\in\mathbb R$.
  \item(b)\ There is a finite set $E$ such that for $z\in\mathbb R\backslash E$
  \begin{align*}
  &\mathcal{U}\left(\overline{u},\underline{v}\right)=\overline{u}''(z)-c\overline{u}'(z)+f(\overline{u}(z))\left(p(\overline{u}(z))-\underline{v}(z)\right)\leq0,\\[0.2cm]  &\mathcal{U}\left(\underline{u},\overline{v}\right)=\underline{u}''(z)-c\underline{u}'(z)+f(\underline{u}(z))\left(p(\underline{u}(z))-\overline{v}(z)\right)\geq0,\\[0.2cm] &\mathcal{V}\left(\overline{u},\overline{v}\right)=d\overline{v}''(z)-c\overline{v}'(z)+s\overline{v}(z)\left(1-\frac{\overline{v}(z)}{q(\overline{u}(z))}\right)\leq0,\\[0.2cm]
  &\mathcal{V}\left(\underline{u},\underline{v}\right)=d\underline{v}''(z)-c\underline{v}'(z)+s\underline{v}(z)\left(1-\frac{\underline{v}(z)}{q(\underline{u}(z))}\right)\geq0.
  \end{align*}
\end{description}
\end{definition}

To apply Schauder's fixed point theorem, we define the Banach space
\begin{equation*}
  B_\rho\left(\mathbb{R},\mathbb{R}^2\right)=\left\{\phi(z)=(u(z),v(z))\in C\left(\mathbb{R},\mathbb{R}^2\right): |\phi(z)|_\rho<+\infty \right\}
\end{equation*}
with the norm
\begin{equation*}
  |\phi(z)|_\rho=\max\left\{\sup\limits_{z\in\mathbb{R}}|u(z)|e^{-\rho|z|},\ \sup\limits_{z\in\mathbb{R}}|v(z)|e^{-\rho|z|}\right\}
\end{equation*}
for small positive constant $\rho$. Then we will seek traveling waves in wave profile set
\begin{equation*}
  \Sigma=\left\{(u,v)\in C\left(\mathbb{R},\mathbb{R}^2\right): \underline{u}(z)\leq u(z)\leq\overline{u}(z),\ \underline{v}(z)\leq v(z)\leq\overline{v}(z)\right\}.
\end{equation*}
Obviously, $\Sigma$ is bounded, closed and convex in $C\left(\mathbb{R},\mathbb{R}^2\right)$. Furthermore, we define
\begin{align*}
  &F(u,v)=\tau u+f(u)(p(u)-v),\\[0.2cm]
  &G(u,v)=\tau v+sv\left(1-\frac{v}{q(u)}\right),
\end{align*}
where $\tau$ satisfies
\begin{equation*}
  \tau>\max\left\{\max\limits_{u\in[0,1]}\left\{f(1)|p'(u)|+q(1)|f'(u)|\right\},\ s\left(2\frac{q(1)}{q(0)}-1\right)\right\}.
\end{equation*}
At the moment, one can verify that for $(u,v)\in[0,1]\times[0,q(1)]$
\begin{align*}
  \frac{\partial F(u,v)}{\partial u}\geq0,\quad \frac{\partial F(u,v)}{\partial v}\leq0,\quad
  \frac{\partial G(u,v)}{\partial u}\geq0,\quad \frac{\partial G(u,v)}{\partial v}\geq0.
\end{align*}
Thus, system \eqref{eq:tr} can be rewritten as
\begin{equation*}
\left\{\begin{split}
&u''-cu'-\tau u+F(u,v)=0,\\[0.2cm]
&dv''-cv'-\tau v+G(u,v)=0.
\end{split}
\right.
\end{equation*}
Next, we define the operator $P=(P_1,P_2): \Sigma\rightarrow C\left(\mathbb{R},\mathbb{R}^2\right)$ by
\begin{align*}
P_1(u,v)(z)&=\frac{1}{\nu_1^+-\nu_1^-}\left\{\int_{-\infty}^{z}e^{\nu_1^-(z-y)} +\int_{z}^{\infty}e^{\nu_1^+(z-y)}\right\} F(u,v)(y)dy,\\[0.2cm]
P_2(u,v)(z)&=\frac{1}{d\left(\nu_2^+-\nu_2^-\right)}\left\{\int_{-\infty}^{z}e^{\nu_2^-(z-y)} +\int_{z}^{\infty}e^{\nu_2^+(z-y)}\right\} G(u,v)(y)dy,
\end{align*}
where
\begin{align*}
\nu_1^\pm=\frac{1}{2}\left(c\pm\sqrt{c^2+4\tau}\right),\qquad \nu_2^\pm=\frac{1}{2d}\left(c\pm\sqrt{c^2+4d\tau}\right).
\end{align*}
It is easy to check that a fixed point of operator $P$ is identical to the solution of system \eqref{eq:tr}. Then a standard argument similar to \cite{w01,z16,m01,r17} gives that $P(\Sigma)\subseteq\Sigma$, $P$ is continuous and compact with respect to the norm in $B_\rho\left(\mathbb{R},\mathbb{R}^2\right)$. Thus, solution of system \eqref{eq:tr} is found by using  Schauder's fixed point theorem. As a consequence, we conclude the following lemma.

\begin{lemma}\label{le:ssm}
Assume that system \eqref{eq:tr} has a pair of upper and lower solutions $\left(\overline{u}(z),\overline{v}(z)\right)$ and $\left(\underline{u}(z),\underline{v}(z)\right)$ satisfying
\begin{equation*}
    \underline{u}(z_i^-)\leq\underline{u}(z_i^+),\ \overline{u}(z_i^+)\leq\overline{u}(z_i^-), \ \underline{v}(z_i^-)\leq\underline{v}(z_i^+),\ \overline{v}(z_i^+)\leq\overline{v}(z_i^-)\ \text{for}\ z_i\in E.
\end{equation*}
Then it has a solution $(u(z),v(z))$ satisfying for $z\in\mathbb R$
\begin{equation*}
  \underline{u}(z)\leq u(z)\leq\overline{u}(z),\quad \underline{v}(z)\leq v(z)\leq\overline{v}(z).
\end{equation*}
\end{lemma}
\section{The existence of semi-traveling wave}
\noindent

In this section, we employ the upper and lower solutions method to establish a positive solution of system \eqref{eq:tr} that meets the left-hand tail limit in \eqref{eq:b}, that is, semi-traveling wave. We further investigate the asymptotic behavior of such waves at $z=-\infty$ by using the unstable manifold theorem, and then provide the sufficient and necessary condition for the existence of semi-traveling wave.

To construct an appropriate pair of upper and lower solutions, we define $P(c,\lambda)$ by
\begin{equation*}
  P(c,\lambda)=d\lambda^2-c\lambda+s,
\end{equation*}
which gives two roots $\lambda_1$ and $\lambda_2$ as follows
\begin{align}\label{eq:en}
  \lambda_1=\frac{1}{2d}\left(c-\sqrt{c^2-4ds}\right),\quad \lambda_2=\frac{1}{2d}\left(c+\sqrt{c^2-4ds}\right).
\end{align}

\textbf{Case 1:} Assume that $c>c^*=2\sqrt{ds}$, then $0<\lambda_1<\lambda_2$. First, for given constants
\begin{align*}
  &0<\beta<\min\left\{c,\ \lambda_1\right\}, \\[0.2cm]
  &0<\varepsilon<\min\left\{\lambda_1,\ \lambda_2-\lambda_1\right\},\\[0.2cm]
  &\sigma>\max\left\{1,\ \frac{f(1)q(1)}{\beta\left(c-\beta\right)}\right\},\\[0.2cm]
  &r>\max\left\{1,\ \frac{-sq(1)}{P(c,\lambda_1+\varepsilon)q(0)}\right\},
\end{align*}
we define functions
\begin{align}\label{eq:sup-sub1}
\begin{split}
  \overline{u}(z)&=1,
  \qquad\qquad \qquad\quad \ \ \ \
  \underline{u}(z)=\left\{
  \begin{array}{ll}
  1-\sigma e^{\beta z}, &z\leq z_1,\\[0.2cm]
  0, &z>z_1,
  \end{array}
  \right.
  \\[0.4cm]
  \overline{v}(z)&=\left\{
  \begin{array}{ll}
  q(1)e^{\lambda_1z}, &z\leq 0,\\[0.2cm]
  q(1), &z> 0,
  \end{array}
  \right.
  \quad
  \underline{v}(z)=\left\{
  \begin{array}{ll}
  q(1)e^{\lambda_1z}\left(1-re^{\varepsilon z}\right), &z\leq z_2,\\[0.2cm]
  0, &z>z_2,
  \end{array}
  \right.
\end{split}
\end{align}
where $z_1=-(1/\beta)\ln \sigma<0$ and $z_2=-(1/\varepsilon)\ln r<0$ by the choice of $\sigma$ and $r$.\\

\textbf{Case 2:} Assume that $c=c^*$, we have $\lambda_1=\lambda_2=\lambda:=c/\left(2d\right)$. For given constants
\begin{align*}
  &h=\lambda e^2/2,\\[0.2cm]
  &0<\beta<\min\left\{c,\ \lambda\right\},\\[0.2cm]
  &\sigma>\max\left\{e^{2\beta/\lambda},\ \frac{f(1)hq(1)}{\left(c-\beta\right)\left(\lambda-\beta\right)\beta e}\right\},\\[0.2cm]
  &r>\max\left\{h\sqrt{2/\lambda},\ \frac{4sh^2q(1) }{d q(0)}\left(\frac{7}{2e\lambda}\right)^{7/2}\right\},
\end{align*}
we define functions
\begin{align}\label{eq:sup-sub2}
\begin{split}
  \overline{u}(z)&=1,
  \qquad\qquad\qquad\qquad\qquad\quad\ \quad
  \underline{u}(z)=\left\{
  \begin{array}{ll}
  1-\sigma e^{\beta z},  &z\leq z_1\\[0.2cm]
  0, &z>z_1,
  \end{array}
  \right.
  \\[0.4cm]
  \overline{v}(z)&=\left\{
  \begin{array}{ll}
  -hq(1)ze^{\lambda z}, &z\leq -2/\lambda,\\[0.2cm]
  q(1), &z> -2/\lambda,
  \end{array}
  \right.
  \quad
  \underline{v}(z)=\left\{
  \begin{array}{ll}
  q(1)e^{\lambda z}\left(-hz-r(-z)^{1/2}\right), &z\leq z_2,\\[0.2cm]
  0, &z>z_2.
  \end{array}
  \right.
\end{split}
\end{align}
where $z_1=-(1/\beta)\ln \sigma<-2/\lambda$ and $z_2=-(r/h)^2<-2/\lambda$ by the choice of $\sigma$ and $r$.

Next, we will check that $(\overline{u}(z),\overline{v}(z))$ and $(\underline{u}(z),\underline{v}(z))$ described in \eqref{eq:sup-sub1} and \eqref{eq:sup-sub2} are indeed a pair of upper and lower solutions of system \eqref{eq:tr} for $c\geq c^*$. Apparently, $\underline{u}(z)\leq\overline{u}(z)$, $\underline{v}(z)\leq\overline{v}(z)$ for $z\in\mathbb R$, then it is sufficient to verify the requirement \emph{(b)} in Definition \ref{Def}.\\

\textbf{Verification of upper and lower solutions for $c>c^*$:} First, since $\overline{u}\equiv1$, $\underline{v}\geq0$ and $p(1)=0$, we have for $z\in\mathbb{R}$
\begin{equation*}
  \mathcal{U}\left(\overline{u},\underline{v}\right)=-f(1)\underline{v}\leq0.
\end{equation*}
Next, for $z>z_1$, it follows from $\underline{u}=0$ that $\mathcal{U}\left(\underline{u},\overline{v}\right)=0$. For $z<z_1<0$, since $\underline{u}=1-\sigma e^{\beta z}$ and $\overline{v}=q(1)e^{\lambda_1z}$,
we have
\begin{align*}
  \mathcal{U}\left(\underline{u},\overline{v}\right)=& e^{\beta z}\left[\beta\sigma\left(c-\beta\right)-f(1-\sigma e^{\beta z})q(1)e^{(\lambda_1-\beta)z}\right]+f(1-\sigma e^{\beta z})p(1-\sigma e^{\beta z})\\[0.2cm]
  \geq& e^{\beta z}\left[\beta\sigma\left(c-\beta\right)-f(1)q(1)\right]>0
\end{align*}
by the choice of $\beta$ and $\sigma$.
For $\mathcal{V}\left(\overline{u},\overline{v}\right)$, we know that $\overline{v}=q(1)$ for $z>0$, then
\begin{equation*}
  \mathcal{V}\left(\overline{u},\overline{v}\right)=\mathcal{V}\left(1,q(1)\right)=sq(1)\left(1-\frac{q(1)}{q(1)}\right)=0.
\end{equation*}
For $z<0$, $\overline{v}=q(1)e^{\lambda_1z}$. From  the definition of $\lambda_1$, we have
\begin{align*}
  \mathcal{V}\left(\overline{u},\overline{v}\right)&=q(1)d\lambda_1^2e^{\lambda_1z}-q(1)c\lambda_1e^{\lambda_1z}+s\lambda_1e^{\lambda_1z}-sq(1)e^{2\lambda_1z}\\[0.2cm]
  &=P(c,\lambda_1)q(1)e^{\lambda_1z}-sq(1)e^{2\lambda_1z}\\[0.2cm]
  &=-sq(1)e^{2\lambda_1z}\leq0.
\end{align*}
Finally, for $z>z_2$, $\mathcal{V}\left(\underline{u},\underline{v}\right)=0$ since $\underline{v}=0$. Note that $\underline{v}=q(1)e^{\lambda_1z}\left(1-re^{\varepsilon z}\right)$ and $\overline{v}=q(1)e^{\lambda_1z}$ for $z<z_2<0$, we have for $z<z_2$
\begin{align*}
  &\underline{v}=\overline{v}-rq(1)e^{\left(\lambda_1+\varepsilon\right)z},\\[0.2cm]
  &\underline{v}'=\overline{v}'-rq(1) \left(\lambda_1+\varepsilon\right)e^{\left(\lambda_1+\varepsilon\right)z},\\[0.2cm]
  &\underline{v}''=\overline{v}''-rq(1) \left(\lambda_1+\varepsilon\right)^2e^{\left(\lambda_1+\varepsilon\right)z}.
\end{align*}
Then by the choice of $\varepsilon$ and $r$, we further deduce that
\begin{align*}
  \mathcal{V}\left(\underline{u},\underline{v}\right)=&d\overline{v}''-drq(1) \left(\lambda_1+\varepsilon\right)^2e^{\left(\lambda_1+\varepsilon\right)z}-c\overline{v}'+crq(1) \left(\lambda_1+\varepsilon\right)e^{\left(\lambda_1+\varepsilon\right)z}\\[0.2cm]
  &+s\overline{v}-srq(1)e^{\left(\lambda_1+\varepsilon\right)z}-\frac{s(\overline{v}-rq(1)e^{\left(\lambda_1+\varepsilon\right)z})^2}{q(\underline{u})}\\[0.2cm]
  \geq&q(1)e^{\lambda_1z}\left(d\lambda_1^2-c\lambda_1+s\right)-rq(1)e^{\left(\lambda_1+\varepsilon\right)z}\left[d \left(\lambda_1+\varepsilon\right)^2-c\left(\lambda_1+\varepsilon\right)+s\right]-\frac{s\overline{v}^2}{q(0)}\\[0.2cm]
  =&q(1)e^{\lambda_1z}P(c,\lambda_1)-rq(1)e^{\left(\lambda_1+\varepsilon\right)z}P(c,\lambda_1+\varepsilon)-\frac{s\overline{v}^2}{q(0)}\\[0.2cm]
  \geq&-rq(1)e^{\left(\lambda_1+\varepsilon\right)z}P(c,\lambda_1+\varepsilon)-\frac{s\overline{v}^2}{q(0)}\\[0.2cm]
  \geq&q(1)e^{\left(\lambda_1+\varepsilon\right)z}\left[-rP(c,\lambda_1+\varepsilon)-\frac{s q(1)}{q(0)}e^{(\lambda_1-\varepsilon)z}\right]\\[0.2cm]
  \geq&q(1)e^{\left(\lambda_1+\varepsilon\right)z}\left[-rP(c,\lambda_1+\varepsilon)-\frac{s q(1)}{q(0)}\right]>0.
\end{align*}

\textbf{Verification of upper and lower solutions for $c=c^*$:} First, as before, we have
\begin{align*}
  &\mathcal{U}\left(\overline{u},\underline{v}\right)=-f(1)\underline{v}\leq0, \text{ for } z\in\mathbb{R},\\[0.2cm]
  &\mathcal{U}\left(\underline{u},\overline{v}\right)=0, \text{ for } z>z_1,\\[0.2cm]
  &\mathcal{V}\left(\overline{u},\overline{v}\right)=0, \text{ for } z>-2/\lambda,\\[0.2cm]
  &\mathcal{V}\left(\underline{u},\underline{v}\right)=0, \text{ for } z>z_2.
\end{align*}
Let us go on to prove the remaining inequalities. For $\mathcal{U}\left(\underline{u},\overline{v}\right)$, when $z<z_1<-2/\lambda$, we have $\underline{u}=1-\sigma e^{\beta z}$ and $\overline{v}=-hq(1)ze^{\lambda z}$.
Note that for $z\in\mathbb R$
\begin{equation*}
  z e^{(\lambda-\beta)z}\geq-\frac{1}{(\lambda-\beta)e},
\end{equation*}
after a simple computation, we arrive at
\begin{align*}
  \mathcal{U}\left(\underline{u},\overline{v}\right)=&e^{\beta z}\left[\beta\sigma\left(c-\beta\right)+f(1-\sigma e^{\beta z})hq(1)z e^{(\lambda-\beta)z}\right]+f(1-\sigma e^{\beta z})p(1-\sigma e^{\beta z})\\[0.2cm]
  \geq& e^{\beta z}\left[\beta\sigma\left(c-\beta\right)+f(1-\sigma e^{\beta z})hq(1)z e^{(\lambda-\beta)z}\right]\\[0.2cm]
  \geq& e^{\beta z}\left[\beta\sigma\left(c-\beta\right)+f(1)hq(1)z e^{(\lambda-\beta)z}\right]\\[0.2cm]
  \geq& e^{\beta z}\left[\beta\sigma\left(c-\beta\right)-\frac{f(1)hq(1)}{(\lambda-\beta)e}\right]>0
\end{align*}
by the choice of $\beta$ and $\sigma$. Next, for $\mathcal{V}\left(\overline{u},\overline{v}\right)$, it is easy to see that if $z<-2/\lambda$, then
\begin{equation*}
  \overline{v}'=-hq(1)\left(1+\lambda z \right) e^{\lambda z},\quad \overline{v}''=-hq(1)\lambda\left(2+\lambda z \right) e^{\lambda z}.
\end{equation*}
Then we deduce that for $z<-2/\lambda$
\begin{equation}\label{eq:u}
\begin{split}
  \mathcal{V}\left(\overline{u},\overline{v}\right)\leq&d\overline{v}''-c\overline{v}'+s\overline{v}\\[0.2cm]
  =&-dhq(1)\lambda\left(2+\lambda z\right) e^{\lambda z}+chq(1)\left(1+\lambda z\right) e^{\lambda z}-shq(1)z e^{\lambda z}\\[0.2cm]
  =&-hq(1)e^{\lambda z}\left[ P(c,\lambda)z+(2d\lambda-c)\right]=0.
\end{split}
\end{equation}
Finally, for $z<z_2<-2/\lambda$, $\overline{v}=-hq(1)ze^{\lambda z}$, we have
\begin{equation*}
  \underline{v}=q(1)e^{\lambda z}\left(-hz-r(-z)^{1/2}\right)=\overline{v}-rq(1)(-z)^{1/2}e^{\lambda z}.
\end{equation*}
Direct calculation gives
\begin{align*}
  &\underline{v}'=\overline{v}'+rq(1)e^{\lambda z}\left[\frac{1}{2}(-z)^{-1/2}-\lambda(-z)^{1/2}\right],\\[0.2cm]
  &\underline{v}''=\overline{v}''+rq(1)e^{\lambda z}\left[\lambda(-z)^{-1/2}+\frac{1}{4}(-z)^{-3/2}-\lambda^2(-z)^{1/2}\right].
\end{align*}
Substituting the above equalities into $\mathcal{V}\left(\underline{u},\underline{v}\right)$, we obtain
\begin{align*}
  \mathcal{V}\left(\underline{u},\underline{v}\right)=&d\overline{v}''+drq(1)e^{\lambda z}\left[\lambda(-z)^{-1/2}+\frac{1}{4}(-z)^{-3/2}-\lambda^2(-z)^{1/2}\right]\\[0.2cm]
  &-c\overline{v}'-crq(1)e^{\lambda z}\left[\frac{1}{2}(-z)^{-1/2}-\lambda(-z)^{1/2}\right]+s\overline{v}-srq(1)(-z)^{1/2}e^{\lambda z}-\frac{s\underline{v}^2 }{q(\underline{u})}.
\end{align*}
From \eqref{eq:u}, we have for $z<-2/\lambda$
\begin{equation*}
  d\overline{v}''-c\overline{v}'+s\overline{v}=0.
\end{equation*}
Since $q(u)$ is increasing function for $u>0$, then by the choice of $r$, we further deduce that
\begin{align*}
  \mathcal{V}\left(\underline{u},\underline{v}\right)=&drq(1)e^{\lambda z}\left[\lambda(-z)^{-1/2}+\frac{1}{4}(-z)^{-3/2}-\lambda^2(-z)^{1/2}\right]\\[0.2cm]
  &-crq(1)e^{\lambda z}\left[\frac{1}{2}(-z)^{-1/2}-\lambda(-z)^{1/2}\right]-srq(1)(-z)^{1/2}e^{\lambda z}-\frac{s\underline{v}^2 }{q(\underline{u})},\\[0.2cm]
  =&q(1)e^{\lambda z}\left[\frac{1}{4}dr(-z)^{-3/2}+r(-z)^{-1/2}\left(d\lambda-\frac{c}{2}\right)-r(-z)^{1/2}P(c,\lambda)\right]-\frac{s\underline{v}^2 }{q(\underline{u})}\\[0.2cm]
  \geq&q(1)e^{\lambda z}\left[\frac{1}{4}dr(-z)^{-3/2}-\frac{s h^2q(1)}{q(0)}z^2e^{\lambda z}\right]\\[0.2cm]
  =&q(1)e^{\lambda z}(-z)^{-3/2}\left[\frac{1}{4}dr-\frac{s h^2q(1)}{q(0)}(-z)^{7/2}e^{\lambda z}\right]\\[0.2cm]
  \geq&q(1)e^{\lambda z}(-z)^{-3/2}\left[\frac{1}{4}dr-\frac{sh^2q(1) }{q(0)}\left(\frac{7}{2e\lambda}\right)^{7/2}\right]>0,
\end{align*}
by taking advantage of the fact that for $z\leq0$
\begin{equation*}
  (-z)^{7/2}e^{\lambda z}\leq\left(\frac{7}{2e\lambda}\right)^{7/2}.
\end{equation*}

Summarizing the above discussion, we further conclude the following lemma.
\begin{lemma}\label{le:ul}
 For $c\geq c^*$, system \eqref{eq:tr} has a non-negative solution $(u(z),v(z))$.
\end{lemma}
\Proof For $c\geq c^*$, function pairs $(\overline{u}(z),\overline{v}(z))$ and $(\underline{u}(z),\underline{v}(z))$ described in \eqref{eq:sup-sub1} or \eqref{eq:sup-sub2} are a pair of upper and lower solutions of system \eqref{eq:tr}, respectively. Moreover, for $c>c^*$, we have
  \begin{align*}
  &\underline{u}'(z_1^-)=-\beta<0=\underline{u}'(z_1^+),\\[0.2cm]
  &\overline{v}'(0^+)=0<q(1)\lambda_1=\overline{v}'(0^-),\\[0.2cm]
  &\underline{v}'(z_2^-)=-q(1)\varepsilon e^{\lambda_1 z_2}<0=\underline{v}'(z_2^+).
\end{align*}
While, for $c=c^*$, we have
  \begin{align*}
  &\underline{u}'(z_1^-)=-\beta<0=\underline{u}'(z_1^+),\\[0.2cm]
  &\underline{v}'(z_2^-)=-q(1)he^{\lambda z_2}/2<0=\underline{v}'(z_2^+),\\[0.2cm]
  &\overline{v}'((-2/\lambda)^+)=0<q(1)he^{-2}=\overline{v}'((-2/\lambda)^-).
\end{align*}
Hence, we complete the proof with Lemma \ref{le:ssm}.

\qed

Actually, the above non-negative solution is a semi-traveling wave of system \eqref{eq:h11}.

\begin{theorem}\label{th:41}
For $c\geq c^*:=2\sqrt{ds}$, system \eqref{eq:tr} has a solution $(u(z),v(z))$ satisfying
\begin{equation*}
  (u(-\infty),u'(-\infty),v(-\infty),v'(-\infty))=(1,0,0,0)
\end{equation*}
and $0<u(z)<1$, $0<v(z)<q(1)$ for $z\in\mathbb R$.
\end{theorem}
\Proof From Lemma \ref{le:ssm}, we can easily get that for $c\geq c^*$
\begin{equation*}
  0\leq\underline{u}(z)\leq u(z)\leq\overline{u}(z)=1\text{ and }0\leq\underline{v}(z)\leq v(z)\leq\overline{v}(z)\leq q(1) \text{ for } z\in\mathbb R,
\end{equation*}
And it is clear that
\begin{align*}
&1=\lim_{z\rightarrow-\infty}\underline{u}(z)\leq\liminf_{z\rightarrow-\infty}u(z)\leq\limsup_{z\rightarrow-\infty}u(z)\leq\lim_{z\rightarrow-\infty}\overline{u}(z)=1,\\[0.2cm]
&0=\lim_{z\rightarrow-\infty}\underline{v}(z)\leq\liminf_{z\rightarrow-\infty}v(z)\leq\limsup_{z\rightarrow-\infty}v(z)\leq\lim_{z\rightarrow-\infty}\overline{v}(z)=0,
\end{align*}
which implies that $\left(u(-\infty),v(-\infty)\right)=(1,0)$.

We further claim that $0<u(z)<1$, $0<v(z)<q(1)$ for $z\in\mathbb R$. For contradiction, we assume that $u(z_0)=0$ for some $z_0\in \mathbb{R}$, then $u'(z_0)=0$ owing to $u(z)\geq0$ over $\mathbb{R}$. The uniqueness of solution yields $u(z)\equiv0$, which contradicts $u(z)\geq\underline{u}(z)>0$ for $z<z_1$. Hence, we have $u(z)>0$ over $\mathbb{R}$. And similar reasons lead to $v(z)>0$ over $\mathbb{R}$. In order to prove $u(z)<1$ over $\mathbb{R}$, we also assume that $u(z_0)=1$ for some $z_0\in \mathbb{R}$. It follows that $u'(z_0)=0$ and $u''(z_0)\leq0$.
Using the first equation of system \eqref{eq:tr} and the assumption \emph{(H1)}-\emph{(H2)}, we have $0=u''(z_0)-f(u(z_0))v(z_0)\leq-f(u(z_0))v(z_0)<0$. Hence, this contradiction yields $u(z)<1$ over $\mathbb{R}$. Similarly, we also get $v(z)<q(1)$ over $\mathbb{R}$ by using $u(z)<1$ over $\mathbb{R}$.

Finally, we show that $(u'(-\infty),v'(-\infty))=(0,0)$. Indeed, for the first equation of system \eqref{eq:tr}, we utilize the variation constants formula to deduce that
\begin{align*}
  &u'(z)=e^{c(z-\varsigma)}u'(\varsigma)-e^{cz}\int_z^\varsigma e^{-c\tau}\left[f(u(\tau))\left(v(\tau)-p(u(\tau))\right)\right] d\tau,
\end{align*}
thus, it is easy to check that for $z\leq \varsigma$
\begin{align*}
  c|u'(z)|\leq ce^{c(z-\varsigma)}|u'(\varsigma)|+\max\limits_{\tau\leq \varsigma}\Big|f(u(\tau))\left(v(\tau)-p(u(\tau))\right)\Big|.
\end{align*}
Consequently, for fixed $\varsigma$, we have
\begin{equation*}
  \limsup\limits_{z\rightarrow-\infty}|u'(z)|\leq\frac{1}{c}\left(\max\limits_{\tau\leq \varsigma}\Big|f(u(\tau))\left(v(\tau)-p(u(\tau))\right)\Big|\right).
\end{equation*}
It follows from \emph{(H1)} and \emph{(H2)} that
\begin{equation*}
  f(u(-\infty))\left(p(u(-\infty))-v(-\infty)\right)=0,
\end{equation*}
then we have $w(-\infty)=0$ by arbitrary of $\varsigma$. Similar reasons give $v'(-\infty)=0$.

\qed

In order to obtain more detailed information on semi-traveling wave, it is essential for us to study the asymptotic behavior of semi-traveling wave at $z=-\infty$. Hence, we set $w(z)=u'(z)$, $y(z)=v'(z)$ and rewrite system \eqref{eq:tr} as a system of first order ODEs in $\mathbb{R}^4$
\begin{equation}\label{eq:od}
\left\{\begin{split}
  &u'(z)=w(z),\\[0.2cm]
  &w'(z)=cw(z)-f(u(z))(p(u(z))-v(z)),\\[0.2cm]
  &v'(z)=y(z),\\[0.2cm]
  &dy'(z)=cy(z)-sv(z)\left(1-\displaystyle\frac{v(z)}{q(u(z))}\right).
\end{split}
\right.
\end{equation}
The eigenvalues of the linearization of \eqref{eq:od} at $\mathbf{e}_{0}=(1, 0, 0, 0)$ are
\begin{align*}
\lambda_{1} & = \frac{1}{2d}\left(c-\sqrt{c^2-4ds}\right), \quad \lambda_{3}=\frac{1}{2} \left(c+\sqrt{c^{2}-4f(1) p'(1)}\right), \\[0.2cm]
\lambda_{2} & = \frac{1}{2d}\left(c+\sqrt{c^2-4ds}\right), \quad \lambda_{4}=\frac{1}{2} \left(c-\sqrt{c^{2}-4f(1) p'(1)}\right).
\end{align*}
Obviously, for $c\geq c^*$, we have $\lambda_{4}<0<\lambda_{1}, \lambda_{2}, \lambda_{3}$ due to $p'(1)<0<f(1)$.
\begin{lemma}\label{th:42}
For $c\geq c^*$, the positive solution $(u(z),v(z))$ obtained by Theorem \ref{th:41} satisfies
\begin{equation*}
  \lim_{z\rightarrow-\infty} \frac{y(z)}{v(z)}=\lambda_1 \text{ and }
  \lim_{z\rightarrow-\infty} \frac{w(z)}{u(z)-1}=\Lambda\in\left\{\lambda_1, \lambda_3\right\}.
\end{equation*}
\end{lemma}
\Proof Based on the choice of upper and lower solutions, from Lemma \ref{le:ssm}, we have
\begin{equation*}
  v(z)=\left\{
  \begin{array}{ll}
   O(e^{\lambda_1z}),\quad &c>c^*, \\[0.2cm]
   O(-ze^{\lambda_1z}),&c=c^*,
   \end{array}
  \right.
\end{equation*}
as $z\rightarrow-\infty$, which implies that
\begin{equation*}
  \lim_{z\rightarrow-\infty} \frac{y(z)}{v(z)}=\lambda_1.
\end{equation*}
Now, we show the second equality, which is quite complicated and tedious. We only consider the case $c>c^*$, since when $c=c^*$ and $\lambda_{1}=\lambda_{2}<(>,=)\lambda_{3}$, it can be similarly treated by using generalized eigenvectors and the unstable manifold theorem.

(1)\ If $\lambda_{1}, \lambda_{2}\neq \lambda_{3}$, then the corresponding eigenvectors are given by
\begin{equation*}
\mathbf{e}_{1}=\begin{pmatrix}
  -1 \\
  -\lambda_{1} \\
  -\psi\left(\lambda_{1}\right) \\
  -\lambda_{1} \psi\left(\lambda_{1}\right)
\end{pmatrix},
\quad
\mathbf{e}_{2}=\begin{pmatrix}
  -1 \\
  -\lambda_{2} \\
  -\psi\left(\lambda_{2}\right) \\
  -\lambda_{2} \psi\left(\lambda_{2}\right)
\end{pmatrix},
\quad
\mathbf{e}_{3}=\begin{pmatrix}
  -1 \\
  -\lambda_{3} \\
  0 \\
  0
\end{pmatrix},
\end{equation*}
where
\begin{equation*}
  \psi(\lambda)=\frac{\lambda^2-c\lambda+f(1)p'(1)}{f(1)}.
\end{equation*}
Thus, every solution $W(z)$ of the corresponding linearized system, which converges to $\mathbf{e}_{0}$ as $z\rightarrow-\infty$, is given by $$W(z)=\mathbf{e}_{0}+C_1\mathbf{e}_{1}e^{\lambda_1z}+C_2\mathbf{e}_{2}e^{\lambda_2z}+C_3\mathbf{e}_{3}e^{\lambda_3z}$$ for some constants $C_i$, $i=1,2,3$.
Applying the unstable manifold theorem yields that as $z\rightarrow-\infty$, there are $\alpha$, $\beta$ and $\gamma$ such that
\begin{align*}
  u(z)&=1-\alpha e^{\lambda_1z}-\beta e^{\lambda_2z}-\gamma e^{\lambda_3z}+h.o.t.,\\[0.2cm]
  w(z)&=0-\alpha \lambda_1 e^{\lambda_1z}-\beta \lambda_2 e^{\lambda_2z}-\gamma \lambda_3 e^{\lambda_3z}+h.o.t.,\\[0.2cm]
  v(z)&=0-\alpha \psi\left(\lambda_{1}\right) e^{\lambda_1z}-\beta\psi\left(\lambda_{2}\right) e^{\lambda_2z}+h.o.t..
\end{align*}
Note that $v(z)=O(e^{\lambda_1z})$ as $z\rightarrow-\infty$, we must have $\alpha\neq0$. Then the following hold.
\begin{description}
  \item (i) If $\lambda_1<\lambda_3$, then $\psi(\lambda_1)<0$. Hence, $\alpha>0$ due to $v(z)>0$ for $z\in\mathbb R$.
  \item (ii) If $\lambda_1>\lambda_3$, then $\psi(\lambda_1)>0$. Hence, $\alpha<0<\gamma$ due to $v(z)>0$ and $u(z)<1$ for $z\in\mathbb R$.
\end{description}
Consequently, one can easily verify that
\begin{equation*}
  \lim_{z\rightarrow-\infty} \frac{w(z)}{u(z)-1}=\min\left\{\lambda_1, \lambda_3\right\}.
\end{equation*}

(2)\ If $\lambda_{1}=\lambda_{3}<\lambda_{2}$, we have a generalized eigenvector
\begin{equation*}
\mathbf{e}_{4}=\begin{pmatrix}
  0 \\
  -f(1)/(2\lambda_{1}-c) \\
  -1 \\
  -\lambda_{1}
\end{pmatrix}.
\end{equation*}
Thus, every solution $W(z)$ of the corresponding linearized system, which converges to $\mathbf{e}_{0}$ as $z\rightarrow-\infty$, is given by $$W(z)=\mathbf{e}_{0}+C_1\left(\mathbf{e}_{1}z+\mathbf{e}_{4}\right)e^{\lambda_1z}+C_2\mathbf{e}_{2}e^{\lambda_2z}$$ for some constants $C_i$, $i=1,2$. Applying the unstable manifold theorem yields that as $z\rightarrow-\infty$, there are $\alpha$ and $\beta$ such that
\begin{align*}
  u(z)&=1-\alpha ze^{\lambda_1z}-\beta e^{\lambda_2z}+h.o.t.,\\[0.2cm]
  w(z)&=0-\alpha \left[\lambda_1z+f(1)/(2\lambda_{1}-c)\right] e^{\lambda_1z}-\beta \lambda_2 e^{\lambda_2z}+h.o.t.,\\[0.2cm]
  v(z)&=0-\alpha e^{\lambda_1z}-\beta\psi\left(\lambda_{2}\right) e^{\lambda_2z}+h.o.t..
\end{align*}
Note that $v(z)=O(e^{\lambda_1z})$ as $z\rightarrow-\infty$ and $v(z)>0$, we must have $\alpha<0$. Moreover, one can easily verify that
\begin{equation}\label{eq:ff}
  \lim_{z\rightarrow-\infty} \frac{w(z)}{u(z)-1}=\lambda_1.
\end{equation}

(3)\ If $\lambda_{1}<\lambda_{2}=\lambda_{3}$, we have a generalized eigenvector
\begin{equation*}
\mathbf{e}_{5}=\begin{pmatrix}
  0 \\
  -f(1)/(2\lambda_{2}-c) \\
  -1 \\
  -\lambda_{2}
\end{pmatrix}.
\end{equation*}
Thus, every solution $W(z)$ of the corresponding linearized system, which converges to $\mathbf{e}_{0}$ as $z\rightarrow-\infty$, is given by $$W(z)=\mathbf{e}_{0}+C_1\mathbf{e}_{1}e^{\lambda_1z}+C_2(\mathbf{e}_{2}z+\mathbf{e}_{5})e^{\lambda_2z}$$ for some constants $C_i$, $i=1,2$. Applying the unstable manifold theorem yields that as $z\rightarrow-\infty$, there are $\alpha$ and $\beta$ such that
\begin{align*}
  u(z)&=1-\alpha e^{\lambda_1z}-\beta ze^{\lambda_2z}+h.o.t.,\\[0.2cm]
  w(z)&=0-\alpha \lambda_1 e^{\lambda_1z}-\beta \left[\lambda_2 z+f(1)/(2\lambda_2-c)\right]e^{\lambda_2z}+h.o.t.,\\[0.2cm]
  v(z)&=0-\alpha \psi\left(\lambda_{1}\right) e^{\lambda_1z}-\beta e^{\lambda_2z}+h.o.t..
\end{align*}
Similarly, we have $\alpha>0$ and \eqref{eq:ff} still hold. Therefore, we complete the proof.

\qed
\begin{remark}\label{re:r1}
From Theorem \ref{th:41}, we know that $v(z)>0$ and $u(z)<1$ over $\mathbb R$. Combining with Lemma \ref{th:42}, one can infer that $y(z)>0>w(z)$ for $z\ll-1$, which is used to prove the monotonicity of traveling wave connecting $(1,0)$ and $(0,\mu)$ in Theorem \ref{th:no1}.
\end{remark}

Finally, let us end this section by proving Theorem \ref{th:ma}.

\textbf{The proof of  Theorem \ref{th:ma}:} We investigate the linearization of the second equation of system \eqref{eq:tr} at $(1,0)$, which provides two eigenvalues $\lambda_1$ and $\lambda_2$ defined in \eqref{eq:en}. If $|c|<c^*$, then  $\lambda_1$ and $\lambda_2$ are a pair of complex eigenvalues. If $c<-c^*$, then  $\lambda_1<\lambda_2<0$. The former implies that $v(z)$ cannot have the same sign as $z\rightarrow-\infty$, while the latter implies that $v(z)$ is unbounded as $z\rightarrow-\infty$. Since neither is feasible, system \eqref{eq:h11} has no semi-traveling wave if $c<c^*$, which concludes the result by coupling Theorem \ref{th:41}.

\qed
\section{The existence of strong-traveling wave}
\noindent

In this section, our main purpose is to prove the existence of strong-traveling waves of system \eqref{eq:h11}. Firstly, inspired by \cite{hsu12,z16,f15,t19}, we shall develop two derivative estimates to produce an invariant set of system \eqref{eq:od}, which is a required and crucial effort. Then semi-traveling wave obtained in Theorem \ref{th:41} can be viewed as a good candidate for strong-traveling waves. Finally, we deduce the right-hand tail limit of such a solution, which establishes the existence of strong-traveling waves. In the sequence, we always assume that $c\geq c^*$.

We start by recalling LaSalle's Invariance Principle.
\begin{proposition}(\cite{h80})\label{le:lsr}
Consider the following initial value problem
\begin{equation}\label{eq:f}
\left\{\begin{split}
&y'=f(y),\\[0.2cm]
&y(0)=y_0,
\end{split}
\right.
\end{equation}
where $f: \mathbb{R}^n\rightarrow \mathbb{R}^n$ is continuous and satisfies Lipschitz condition. Let $\Sigma\subseteq \mathbb{R}^n$ be an open set. Suppose $y(t,y_0)$ is a solution of equation \eqref{eq:f} which is positively invariant in $\Sigma$. If there is a continuous and bounded below function $V(y): \Sigma \rightarrow \mathbb{R}^n$ such that the orbital derivative of $V(y)$ along $y(t,y_0)$ is non-positive, i.e.,
\begin{equation*}
\frac{dV}{dt}=grad\ V(y) \cdot f (y)\leq0,
\end{equation*}
then the $\omega$-limit set of $y(t,y_0)$ is contained in the largest invariant set $$I=\left\{y\in \Sigma: grad\ V(y) \cdot f (y) = 0\right\}.$$
\end{proposition}

\subsection{Derivative estimates}
In this subsection, we shall develop two derivative estimates of $u(z)$ and $v(z)$.
\begin{lemma}\label{le:ud}
There is a large enough $K>q(1)/c$ such that for $z\in\mathbb{R}$,
\begin{equation*}
  -Kf(u(z))< w(z)< Kf(u(z)).
\end{equation*}
\end{lemma}
\Proof To prove the above inequalities, we define functions
\begin{align*}
   \phi_1(z)=w(z)-Kf(u(z)),\\[0.2cm]
   \phi_2(z)=w(z)+Kf(u(z)).
\end{align*}
It is sufficient to show that $\phi_1(z)<0$ and $\phi_2(z)>0$ over $\mathbb R$. By recalling Theorem \ref{th:41}, we have $0<u(z)<1$ and $\left(u(-\infty),w(-\infty)\right)=(1,0)$. Thus, there are  $K>0$ and $z_0\ll-1$ such that
\begin{align*}
 \phi_1(z)<0\ \ \text{and} \ \ \phi_2(z)>0\ \ \text{for}\ \ z\leq z_0.
\end{align*}

(1)\ We show that $\phi_1(z)<0$ for $z>z_0$ by contradiction. Assume that there is a $z_1>z_0$ such that $\phi_1(z_1)=0$ and $\phi'_1(z_1)\geq0$. Since $u(z)$ and $v(z)$ are bounded, and $w(z_1)=Kf(u(z_1))$, then for large enough $K>0$
\begin{align*}
  0&\leq \phi'_1(z_1)
  =w'(z_1)-Kf'(u(z_1))w(z_1)\\[0.2cm]
  &=cw(z_1)-f(u(z_1))p(u(z_1))+f(u(z_1))v(z_1)-Kf'(u(z_1))w(z_1)\\[0.2cm]
  &=\left[cK-p(u(z_1))+v(z_1)-K^2f'(u(z_1))\right]f(u(z_1))<0.
\end{align*}
Hence, the contradiction yields $\phi_1(z)<0$ for $z>z_0$.

(2)\ We show that $\phi_2(z)>0$ for $z>z_0$ by contradiction.
Assume that there is a $z_1>z_0$ such that $\phi_2(z_1)=0$ and $\phi'_2(z_1)\leq0$.
We claim that there is a $z_2>z_1$ such that $\phi_2(z_2)=0$ and $\phi'_2(z_2)\geq0$. To see this, we further assume that  $\phi_2(z)<0$ for all $z>z_1$, that is, $w(z)<-Kf(u(z))$ for all $z>z_1$.
Note that $0<u(z)<1$, $0<v(z)<q(1)$ over $\mathbb R$ and $p(u)>0$ for $u\in(0,1)$, then we have the estimate for $z>z_1$
\begin{align*}
  w'(z)=&cw(z)+f(u(z))v(z)-f(u(z))p(u(z))\\[0.2cm]
  <&[-cK+v(z)-p(u(z))]f(u(z))\\[0.2cm]
  <&(q(1)-cK)f(u(z))<0
\end{align*}
as long as $K>q(1)/c$. Thus, it is easy to check that for $z>z_1$
\begin{equation*}
  w(z)<w(z_1)=-Kf(u(z_1))<0,
\end{equation*}
which contradicts the boundedness of $u(z)$. Hence, the claim is valid. On the other hand, at $z=z_2$, $w(z_2)=-Kf(u(z_2))$, we also get for large enough $K>0$
\begin{align*}		
0&\leq \phi'_2(z_2)=w'(z_2)+Kf'(u(z_2))w(z_2)\\[0.2cm]
&=cw(z_2)-f(u(z_2))p(u(z_2))+f(u(z_2))v(z_2)+Kf'(u(z_2))w(z_2)\\[0.2cm]
&=\left[-(K^2f'(u(z_2))+cK)+v(z_2)-p(u(z_2))\right]f(u(z_2))<0,
\end{align*}
due to the boundedness of $u(z)$ and $v(z)$. Hence, the contradiction yields $\phi_2(z)>0$ for $z>z_0$.

\qed

\begin{lemma}\label{le:vd}
If $c\geq c^*=2\sqrt{ds}$, then for $z\in\mathbb{R}$,
\begin{equation*}
  -\frac{sq(1)}{cq(0)}v(z)\leq y(z)\leq\frac{c}{2d}v(z).
\end{equation*}
\end{lemma}
\Proof To prove the above inequalities, we define function
\begin{align*}
  &\psi_1(z)=y(z)-\frac{c}{2d}v(z),\\[0.2cm]
  &\psi_2(z)=y(z)+\frac{sq(1)}{cq(0)}v(z).
\end{align*}
It is sufficient to show that $\psi_1(z)\leq0$ and $\psi_2(z)\geq0$ over $\mathbb R$.

(1)\ We show that $\psi_1(z)\leq0$ over $\mathbb R$ by contradiction. Assume that there is a $z_0$ such that $\psi_1(z_0)>0$, then the $v(z)$ equation of system \eqref{eq:tr} gives for $c\geq c^*=2\sqrt{ds}$
\begin{align*}
  d\psi_1'(z)-\frac{c}{2}\psi_1(z)&=dy'(z)-cy(z)+\frac{c^2}{4d}v(z)\\[0.2cm]
  &=\left(\frac{sv(z)}{q(u(z))}+\frac{c^2}{4d}-s\right)v(z)\\[0.2cm]
  &\geq\frac{sv^2(z)}{q(u(z))}>0.
\end{align*}
Hence, by the comparison principle, we have $\psi_1(z)>e^{c(z-z_0)/2d}\psi_1(z_0)>0$ for $z>z_0$. Using the expression of $\psi_1(z)$ and the comparison principle again, we obtain
\begin{equation*}
  v(z)>e^{c(z-z_0)/2d}v(z_0) \ \text{for}\ z>z_0.
\end{equation*}
Hence, we have $v(z)\rightarrow +\infty$ as $z\rightarrow+\infty$, which contradicts the boundedness of $v(z)$.

(2)\ We show that $\psi_2(z)\geq0$ over $\mathbb R$. From Lemma \ref{th:42}, we have for $c\geq c^*$
\begin{equation*}
  \lim_{z\rightarrow-\infty} \frac{y(z)}{v(z)}=\lambda_1>0,
\end{equation*}
then there is a $z_0\ll-1$ such that for $z\leq z_0$
\begin{equation*}
  \psi_2(z)=v(z)\left(\frac{y(z)}{v(z)}+\frac{sq(1)}{cq(0)}\right)>0.
\end{equation*}
For contradiction, we assume that there is a $z_1>z_0$ such that $\psi_2(z_1)=0$ and $\psi'_2(z_1)\leq0$. Similar to the proof of Lemma \ref{le:ud},
$\psi_2(z_2)=0$ and $\psi'_2(z_2)\geq0$ for some $z_2>z_1$. Actually, since $0<u(z)<1$, $0<v(z)<q(1)$ over $\mathbb R$ and $q(u(z))>q(0)$, then as long as $y'(z)\leq -sq(1)v(z)/cq(0)$, we have
\begin{equation}\label{eq:ss}
\begin{split}
  dy'(z)&=cy(z)-sv(z)\left(1-\frac{v(z)}{q(u(z))}\right)\\[0.2cm]
  &\leq sv(z)\left(\frac{v(z)}{q(u(z))}-\frac{q(1)}{q(0)}-1\right)<0.
\end{split}
\end{equation}
Therefore, if $\psi_2(z)<0$ for all $z>z_1$, that is, for $z>z_1$
\begin{equation*}
  y(z)<-\frac{sq(1)}{cq(0)}v(z),
\end{equation*}
then \eqref{eq:ss} yields that for $z>z_1$,
\begin{equation*}
  y(z)<y(z_1)=-\frac{sq(1)}{cq(0)}v(z_1)<0,
\end{equation*}
which contradicts the positivity of $v(z)$. Thanks to $\psi'_2(z_2)\geq0$ and $y(z_2)=-sq(1)v(z_2)/cq(0)$, we have
\begin{equation*}
  y'(z_2)\geq -\frac{sq(1)}{cq(0)}y(z_2)=\left(\frac{sq(1)}{cq(0)}\right)^2v(z_2)>0.
\end{equation*}
Thus, from $0<v(z_2)<q(1)$ and $q(u(z_2))>q(0)$, we deduce that
\begin{align*}
  0&=dy'(z_2)-cy(z_2)+sv(z_2)\left(1-\frac{v(z_2)}{q(u(z_2))}\right)\\[0.2cm]
  &\geq sv(z_2)\left(\frac{q(1)}{q(0)}-\frac{v(z_2)}{q(u(z_2))}+1\right)>0,
\end{align*}
which leads to a contradiction. Hence, we complete the proof.

\qed

Then according to Lemmas \ref{le:ud} and \ref{le:vd}, the set
\begin{equation*}
  D=\left\{\chi(z)\Big|
    0<u<1,\ 0<v<q(1), \displaystyle-Kf(u)< w< Kf(u),\
    \displaystyle-\frac{2sq(1)}{cq(0)}v<y<\frac{c}{d}v
  \right\}
\end{equation*}
is an invariant set of system \eqref{eq:od}, where $\chi(z):=(u,w,v,y)(z)$.
\subsection{Convergence of semi-traveling wave}
\noindent

In this subsection, we will deduce the right-hand tail limit of semi-traveling wave obtained in Theorem \ref{th:41}. To be specific, based on the open set $D$, the appropriate Lyapunov functions are constructed for system \eqref{eq:od},  such that the candidate is a strong-traveling wave.

We firstly prove Theorem \ref{th:ma2}. For convenience, let us recall it again.

\begin{theorem}\label{th:ex}
Assume that (H1)-(H2) hold and $h(u)=1$. For $c\geq c^*$, if  $p(0)>\mu$ and
\begin{description}
  \item (P) $(u^*+\mu)f(u)-(u-u^*)(u+\mu)f'(u)>0$ for $u\in(0,1)$,
\end{description}
then system \eqref{eq:h11} has a traveling wave connecting $(1,0)$ and $(u^*,v^*)$.
\end{theorem}
\Proof
 Let
\begin{equation*}
  H(u,v)=\int_{u^*}^u\frac{\eta-u^*}{q(\eta)f(\eta)}d\eta+\frac{1}{s}\int_{v^*}^v\frac{\eta-v^*}{\eta}d\eta.
\end{equation*}
Obviously, $H(u,v)\geq0$ for $u, v>0$.
Then we define Lyapunov function $\mathcal{L}:\mathbb{R}^4\rightarrow\mathbb{R}$ by
\begin{equation*}
  \mathcal{L}(u,w,v,y)=cH(u,v)-w\partial_u H-dy\partial_v H,
\end{equation*}
thus, $\mathcal{L}(u,w,v,y)$ is continuous function with a lower bound since for $(u,w,v,y)\in D$
\begin{align*}
 &w\partial_u H=\frac{w(u-u^*)}{f(u)q(u)}\leq\frac{w(u+u^*)}{q(0)f(u)}\leq\frac{2w}{q(0)f(u)}<\frac{2K}{q(0)},\\[0.2cm]
 &y\partial_v H=\frac{y(v-v^*)}{sv}\leq\frac{1}{s}\left(\frac{c}{d}v-\frac{y}{v}v^*\right)\leq \frac{q(1)}{s}\left(\frac{c}{d}+\frac{2sv^*}{cq(0)}\right).
\end{align*}
From \eqref{eq:od}, the orbital derivative of $\mathcal{L}$ along $\chi(z)$ is
\begin{equation*}
\begin{split}
  \mathcal{L}'(\chi(z))&=(cw-w')\partial_u H+(cy-dy')\partial_v H-w^2 \partial_{uu} H-dy^2 \partial_{vv} H\\[0.2cm]
  &=f(u)(p(u)-v)\partial_u H+sv\left(1-v/q(u)\right)\partial_v H-w^2 \partial_{uu} H-dy^2 \partial_{vv} H\\[0.2cm]
  &=\left[\left(u-u^*\right)(p(u)-v)+\left(v-v^*\right)\left(q(u)-v\right)\right]/q(u)-w^2 \partial_{uu} H-dy^2 \partial_{vv} H.
\end{split}
\end{equation*}
Using the expression of $p(u^*)=q(u^*)=v^*$, we further obtain
\begin{align*}
  &\left(u-u^*\right)(p(u)-v)=\left(u-u^*\right)(p(u)-p(u^*))+\left(u-u^*\right)(v^*-v),\\[0.2cm]
  &\left(v-v^*\right)\left(q(u)-v\right)=\left(v-v^*\right)\left(q(u)-q(u^*)\right)-\left(v-v^*\right)^2.
\end{align*}
Hence,
\begin{align*}
\mathcal{L}'(\chi(z))=&\left[\left(u-u^*\right)(p(u)-p(u^*))-\left(v-v^*\right)^2\right]/q(u)
-w^2 \partial_{uu} H-dy^2 \partial_{vv} H+G(u,v),
\end{align*}
where
\begin{align*}
  G(u,v)&=(v-v^*)\left[(q(u)-q(u^*))-(u-u^*)\right]/q(u)\\[0.2cm]
  &=(v-v^*)\left[(u-u^*)-(u-u^*)\right]/q(u)=0,
\end{align*}
due to $q(u)=uh(u)+\mu=u+\mu$. Moreover, a simple computation gives
\begin{equation*}
\partial_{vv} H=\frac{v^*}{sv^2}\ \text{  and  }\ \partial_{uu} H=\frac{(u^*+\mu)f(u)-(u-u^*)(u+\mu)f'(u)}{(u+\mu)^2f^2(u)}.
\end{equation*}
Thus, \emph{(H2)} and \emph{(P)} imply that the first and third terms are non-positive in $\mathcal{L}'(\chi(z))$.  Hence, the orbital derivative of $\mathcal{L}$ along $\chi(z)$ is non-positive, and it is zero if and only if $\chi(z)=(u^*,0,v^*,0)$. On the other hand, $p(0)>\mu$ implies that there is a unique positive equilibrium $(u^*,v^*)$ in $D$. Therefore, we complete the proof by using LaSalle's Invariance Principle.

\qed

Finally, we show Theorem \ref{th:ma3}. By recalling that
\begin{align*}
  \mathcal{A}=&\left\{v\in C^2\left(\mathbb{R}\right)\mid 0<v(z)<\mu \text{ and } y(z)>0 \text{ over } \mathbb{R}\right\}, \\[0.2cm]
  \mathcal{B}=&\left\{
           v\in C^2\left(\mathbb{R}\right) \bigg|
           \begin{array}{ll}
              \exists \ z_v  \text{ such that } v(z_v)=\mu, y(z)>0 \text{ in } (-\infty, z_v) \\
              \text{and } v(z)>\mu \text{ in } (z_v,+\infty)
           \end{array}
           \right\},
\end{align*}
then Theorem \ref{th:ma3} is repeated as below.

\begin{theorem}\label{th:no1}
Assume that (H1)-(H3) hold.  For $c\geq c^*$, if
\begin{equation*}
  f(u)=ug(u) \quad \text{and} \quad p(u)=(1-u)/g(u),
\end{equation*}
where $g(u)$ satisfies
  \begin{equation}\label{eq:k}
    \min_{u\in[0,1]}g(u)>\frac{h(1)(h(1)+\mu)}{\mu^2h(0)},
  \end{equation}
  then the following hold.
  \begin{description}
    \item (a) System \eqref{eq:h11} has a traveling wave  $(u(z),v(z))$ connecting $(1,0)$ and $(0,\mu)$.
    \item (b) $v(z)\in \mathcal{A}\cup \mathcal{B}$ and $w(z)<0$ over $\mathbb{R}$.
  \end{description}
\end{theorem}
\Proof Firstly, we show the existence by LaSalle's Invariance Principle. It follows from \eqref{eq:k} that there is a $\varrho$ satisfying
\begin{equation}\label{eq:bb}
  \frac{sh(1)}{\mu\left(\min_{u\in[0,1]}g(u)\right)}<\varrho<\frac{s\mu h(0)}{h(1)+\mu}.
\end{equation}
Then we define Lyapunov function $\mathcal{L}:\mathbb{R}^4\rightarrow\mathbb{R}$ by
\begin{align*}
  \mathcal{L}(\chi(z))
  &=\varrho(cu-w)+c\int_{q(0)}^v\frac{\eta-q(0)}{\eta}d\eta-dy+\frac{dq(0) y}{v}.
\end{align*}
Since $q(u)-q(0)=q(u)-\mu=uh(u)$ and \eqref{eq:od}, after a direct computation, we have
\begin{align*}
  \mathcal{L}'(\chi(z))=&\varrho\left(cw-w'\right)+\left(1-q(0)/v\right)(cy-dy')-dq(0) \left(y/v\right)^2\\[0.2cm]
  =&\varrho u\left(1-u-g(u)v\right)+s(v-q(0))\left(1-v/q(u)\right)-dq(0) \left(y/v\right)^2\\[0.2cm]
  =&-\varrho u^2+\varrho u\left(1-g(u)v\right)-s(v-q(0))^2/q(u)\\[0.2cm]
  &+s(v-q(0))(q(u)-q(0))/q(u)-dq(0) \left(y/v\right)^2\\[0.2cm]
  =&-\varrho u^2-s(v-q(0))^2/q(u)-dq(0)\left(y/v\right)^2\\[0.2cm]
  &+\varrho u\left[1-g(u)v+s(v-q(0))h(u)/\varrho q(u)\right]\\[0.2cm]
  =&-\varrho u^2-s(v-q(0))^2/q(u)-dq(0) \left(y/v\right)^2+\varrho uC(u),
\end{align*}
where
\begin{align*}
  C(u)&=1-g(u)v+s(v-q(0))h(u)/\varrho q(u).
\end{align*}
Through the monotonicity of $h(u)$ and $0<u<1$, we have
\begin{align*}
  C(u)
  &=\left[(\varrho q(u)-sq(0)h(u))+v(sh(u)-\varrho q(u)g(u))\right]/\varrho q(u)\\[0.2cm]
  &<\left[(\varrho q(1)-sq(0)h(0))+v(sh(1)-\varrho q(0)g(u))\right]/\varrho q(u)\\[0.2cm]
  &=\left\{\left[\varrho (h(1)+\mu)-s\mu h(0)\right]+v(sh(1)-\varrho \mu g(u))\right\}/\varrho q(u).
\end{align*}
Automatically, $C(u)<0$ for $u\in[0,1]$ by the choice of $\varrho$. Then the orbital derivative of $\mathcal{L}$ along $\chi(z)$ is non-positive and it is zero if and only if $u=0$, $v=\mu$, $y=0$.  Thus, LaSalle's Invariance Principle yields $\chi(z)\rightarrow(0,w(+\infty),\mu,0)$ as $z\rightarrow+\infty$ if $w(+\infty)$ exists. Since $u(z)$, $v(z)$ and $w(z)$ are bounded for $z\in\mathbb R$, then the $u(z)$ equation gives that $w(z)$ is uniformly continuous. And we further infer that $w(+\infty)=0$ by the Barbalat Lemma. Hence, the existence of traveling wave is proved.

Next, we show $v(z)\in \mathcal{A}\cup \mathcal{B}$. We claim that $y(z)\neq0$ as long as $z$ satisfies $0<v(z)\leq \mu$. Actually, we assume that there is a $z_0\in\mathbb R$ such that $y(z_0)=0$ and $0<v(z_0)\leq \mu$. Since $0<u(z_0)<1$, $\mu=q(0)<q(u(z_0))$ by the monotonicity of $q(u)$, we have
\begin{equation*}
  dy'(z_0)=-sv(z_0)\left(1-\frac{v(z_0)}{q(u(z_0))}\right)<-sv(z_0)\left(1-\frac{\mu}{q(0)}\right)=0,
\end{equation*}
which gives that $v(z)$ attains a local maximum at $z_0$. Thus, $0<v(z_1)<\mu$ and $y(z_1)<0$ for some $z_1>z_0$. Observe that
\begin{equation*}
  dy'(z)=cy(z)-sv(z)\left(1-\frac{v(z)}{q(u(z))}\right)<0
\end{equation*}
as long as $z$ satisfies $y(z)\leq0$ and $0<v(z)\leq \mu$, we have $y'(z)<0$ and $y(z)<y(z_1)<0$ for $z>z_1$, which contradicts the positive of $v(z)$. Hence, we prove the assertion of the claim.
Since $v(-\infty)=0$ and $v(+\infty)=\mu$, $y(z)>0$ as long as $z$ satisfies $0<v(z)\leq \mu$. Meanwhile, it holds that either $v(z)\in(0,\mu)$ over $\mathbb{R}$ or there exists a $z_v$ such that $v(z_v)=\mu$.  For the former case, clearly,  $v(z)\in \mathcal{A}$. For the latter case, $y(z)>0$ for $z\leq z_v$. On the other hand, there is no point $z>z_v$ with $v(z)=\mu$ and $y(z)\leq0$, which implies that $v(z)\in \mathcal{B}$.

Finally, we show that $w(z)<0$ over $\mathbb{R}$, which be considered in two different cases.

For the case $v\in\mathcal{A}$. For contradiction, we assume that there is $z_0\in\mathbb R$ such that $w(z_0)=0$.
From Remark \ref{re:r1}, we know that $w(z)<0$ for $z\ll-1$,
then we can reasonably set
\begin{equation*}
  \bar{z}=\sup\left\{z\in(-\infty,z_0]: w(\varsigma)<0 \ \text{for}\ \varsigma\in(-\infty, z)\right\}.
\end{equation*}
It follows from the definition of $\bar{z}$ that $w(\bar{z})=0$ and $w'(\bar{z})\geq0$. Thus, the $u(z)$ equation of system \eqref{eq:tr} gives
\begin{equation}\label{eq:m}
  v(\bar{z})\geq p(u(\bar{z})).
\end{equation}
Using $u(+\infty)=0$ and the definition of $\bar{z}$, we can define
\begin{equation*}
  \hat{z}=\inf\left\{z\in(\bar{z},+\infty): w(\varsigma)>0 \ \text{for}\ \varsigma\in(\bar{z}, z)\right\}.
\end{equation*}
Obviously,
$u(\hat{z})>u(\bar{z})$, $w(\hat{z})=0$ and $w'(\hat{z})\leq0$. Then the $u(z)$ equation gives
\begin{equation}\label{eq:mm}
  v(\hat{z})\leq p(u(\hat{z})).
\end{equation}
Since $p'(u)<0$ for $u>0$, we see that
\begin{equation*}
  v(\bar{z})\geq p(u(\bar{z}))> p(u(\hat{z}))\geq v(\hat{z}),
\end{equation*}
which contradicts with $y(z)>0$ over $\mathbb{R}$. Therefore, we have $w(z)<0$ over $\mathbb{R}$ if $v\in\mathcal{A}$.

For the another case $v\in\mathcal{B}$, we first claim that $w(z)<0$ if $z\geq z_v$. Otherwise, we assume that $w(z_0)\geq0$ for some $z_0\geq z_v$. Since $v(z)\geq \mu$ for $z\geq z_v$, we deduce that for $z\geq z_v$
\begin{equation}\label{eq:vv}
\begin{split}
  w'(z)-cw(z)&=f(u(z))\left(v(z)-p(u(z))\right)\\[0.2cm]
  &\geq f(u(z))\left(\mu-p(0)\right)>0,
\end{split}
\end{equation}
here we have used the fact that
\begin{equation*}
  p(0)=\frac{1}{g(0)}<\frac{1}{\min\limits_{u\in[0,1]}g(u)}<\frac{\mu^2h(0)}{h(1)(h(1)+\mu)}<\mu.
\end{equation*}
For the case $w(z_0)>0$, \eqref{eq:vv} and the comparison principle give for all $z>z_0$
\begin{equation}\label{eq:g}
  w(z)>w(z_0)e^{c\left(z-z_0\right)},
\end{equation}
which contradicts the boundedness of $u(z)$. For the case $w(z_0)=0$, \eqref{eq:g} yields that $w(z_1)>0$ for some $z_1>z_0$, using \eqref{eq:vv} and the comparison principle again, we obtain
all $z>z_1$
\begin{equation*}
  w(z)>w(z_1)e^{c\left(z-z_1\right)},
\end{equation*}
which contradicts the boundedness of $u(z)$.
Thus, we prove the assertion of the claim. Therefore, it is sufficient to prove that $w(z)<0$ for $z<z_v$. Similar to the case $v(z)\in\mathcal{A}$, we assume that there is $z_0<z_v$ such that $w(z_0)=0$ and define
$$\bar{z}=\sup\left\{z\in(-\infty, z_0]: w(\varsigma)<0 \ \text{for}\ \varsigma\in(-\infty, z)\right\}.$$
Obviously, from the definition of $\bar{z}$, we have $w(\bar{z})=0$ and $w'(\bar{z})\geq0$. Then inequality \eqref{eq:m} is still valid. Note that $w(z_v)<0$ from the above claim, then there is a $\hat{z}\in(\bar{z},z_v)$ such that $u(z)$  attains a local maximum at $\hat{z}$, and inequality \eqref{eq:mm} also remains valid. These again contradict the monotonicity of $v(z)$.
Above all, we complete the proof.

\qed

\section{Application}
\noindent

In this section, we apply Theorems \ref{th:ma2} and \ref{th:ma3} to modified Leslie-Gower system with Holling II or Lotka-Volterra type, then some numerical simulations are carried by using MATLAB. It should be pointed out that we achieve a better result on existence of traveling wave connecting $(1,0)$ and $(u^*,v^*)$ by constructing a novel Lyapunov function for system with Lotka-Volterra type.

\subsection{Modified Leslie-Gower system with Holling II type}
\noindent

Let us now consider system
\begin{equation}\label{eq:h2}
\left\{\begin{split}
&\displaystyle\frac{\partial u(x, t)}{\partial t}=\Delta u(x, t)+u(x, t)\left(1-u(x, t)-\frac{av(x, t)}{u(x, t)+e_1}\right), \\[0.2cm]
&\displaystyle\frac{\partial v(x, t)}{\partial t}=d \Delta v(x, t)+sv(x, t)\left(1-\frac{v(x, t)}{u(x, t)+e_2}\right).
\end{split}
\right.
\end{equation}
where $a$, $e_1$, $d$, $s$ and $e_2$ are positive constants. Then
\begin{equation*}
  f(u)=\frac{au}{u+e_1},  \quad p(u)=\frac{1}{a}(1-u)(u+e_1),\quad q(u)=u+e_2, \quad h(u)=1,\quad \mu=e_2.
\end{equation*}

For convenience, we define
\begin{equation*}
  \gamma(x)=(1-x)(x+e_1)/(x+e_2), \ x\in(0,1).
\end{equation*}
By Theorem \ref{th:ma2}, we have the following result on existence of traveling wave connecting $(1,0)$ and $(u^*,v^*)$.
\begin{theorem}\label{th:gg}
For $c\geq c^*$, there is a constant
\begin{align*}
\begin{split}
\bar{a}=
\left\{
\begin{array}{ll}
e_{1}/e_{2}, &\text { if } e_{1} \in[1,+\infty) \cap\left(0, e_{2}\right], \\[0.2cm]
\min \left\{e_{1}/e_{2}, \gamma\left(\displaystyle\frac{e_1-e_2}{1+2e_1+e_1e_2}\right)\right\}, &\text { if } e_{1} \in[1,+\infty) \cap\left(e_{2}, +\infty\right),
\end{array}
\right.
\end{split}
\end{align*}
such that if $a<\bar{a}$, then system \eqref{eq:h2} has a traveling wave connecting $(1,0)$ and $(u^*,v^*)$.
\end{theorem}
\Proof Clearly, Assumption \ref{as} hold if $e_1\geq1$. If $a<e_1/e_2$, there is a $u^*\in(0,1)$ satisfying $K(u)=u^2-(1-a-e_1)u+(ae_2-e_1)=0$, such that system has a unique positive equilibrium $(u^*,v^*)$.
Then we only check the assumption \emph{(P)}, that is, for $u\in(0,1)$
\begin{equation*}
 B(u)=(u^*+e_2-e_1)u^2+2u^*e_1u+e_1e_2u^*>0.
\end{equation*}
\begin{description}
  \item(1) When $u^*=e_1-e_2\in(0,1)$, that is, $e_1\in(e_2,e_2+1)$ and $a=\gamma(e_1-e_2)$. It is obvious that $B(u)>0$ for $u\in(0,1)$.

  \item(2) When $u^*\in(e_1-e_2,1)$, we have $B'(u)>0$ for $u\in(0,1)$ and $B(u)>B(0)=e_1e_2u^*>0$.\\
(a) If $e_1\in(0,e_2]$, we have $u^*\in(0,1)\subseteq(e_1-e_2,1)$.\\
(b) If $e_1\in(e_2,e_2+1)$,  then $u^*\in(e_1-e_2,1)$ is equivalent to $K(e_1-e_2)<0$ from the graph of $K(u)$, that is, $a<\gamma(e_1-e_2)$.
  \item(3) When $u^*\in(0,e_1-e_2)$, $B(u)$ is a concave function. Due to $B(0)>0$, then it is sufficient to prove that $B(1)\geqslant0$, that is
\begin{equation*}
  e_1-e_2>u^*\geqslant\Lambda:=\frac{e_1-e_2}{1+2e_1+e_1e_2}.
\end{equation*}
(a) If $e_1\in(e_2,e_2+1)$, then $u^*\in\left[\Lambda, e_2-e_1\right)$,
 which is equivalent to $K(\Lambda)\leqslant0$ and $K(e_1-e_2)>0$ from the graph of $K(u)$. Note that $\gamma(x)$ is decreasing on $x$ if $e_1>e_2$, thus, we have $\gamma(e_1-e_2)<a\leqslant\gamma(\Lambda)$. \\
(b) If $e_1\in[e_2+1,+\infty)$, then $u^*\in\left[\Lambda, 1\right)$, that is, $a\leqslant\gamma(\Lambda)$.
\end{description}
Based on (1), (2)-(b) and (3), we know that $B(u)>0$ for $u\in(0,1)$ if
\begin{equation*}
  a<\min\left\{e_{1}/e_{2}, \gamma\left(\frac{e_1-e_2}{1+2e_1+e_1e_2}\right)\right\}\ \text{for} \ e_1\in(e_2,+\infty).
\end{equation*}
To summarize, we complete the proof.

\qed

Next, according to Theorem \ref{th:ma3}, we have the following result.
\begin{theorem}\label{th:no3}
Assume that $c\geq c^*$, if $e_1\geq1$ and $ae^2_2>(1+e_1)(1+e_2)$, then system \eqref{eq:h2} has a traveling wave connecting $(1,0)$ and $(0,e_2)$ and satisfying $v(z)\in \mathcal{A}\cup \mathcal{B}$ and $w(z)<0$ over $\mathbb{R}$.
\end{theorem}
\Proof Clearly, if $e_1\geq1$, then Assumption \ref{as} hold. And $g(u)=a/(u+e_1)$. From Theorem \ref{th:ma3}, we complete the proof.

\qed

\begin{figure}[htbp]
  \centering
  \includegraphics[width=2.9in,height=1.7in,clip]{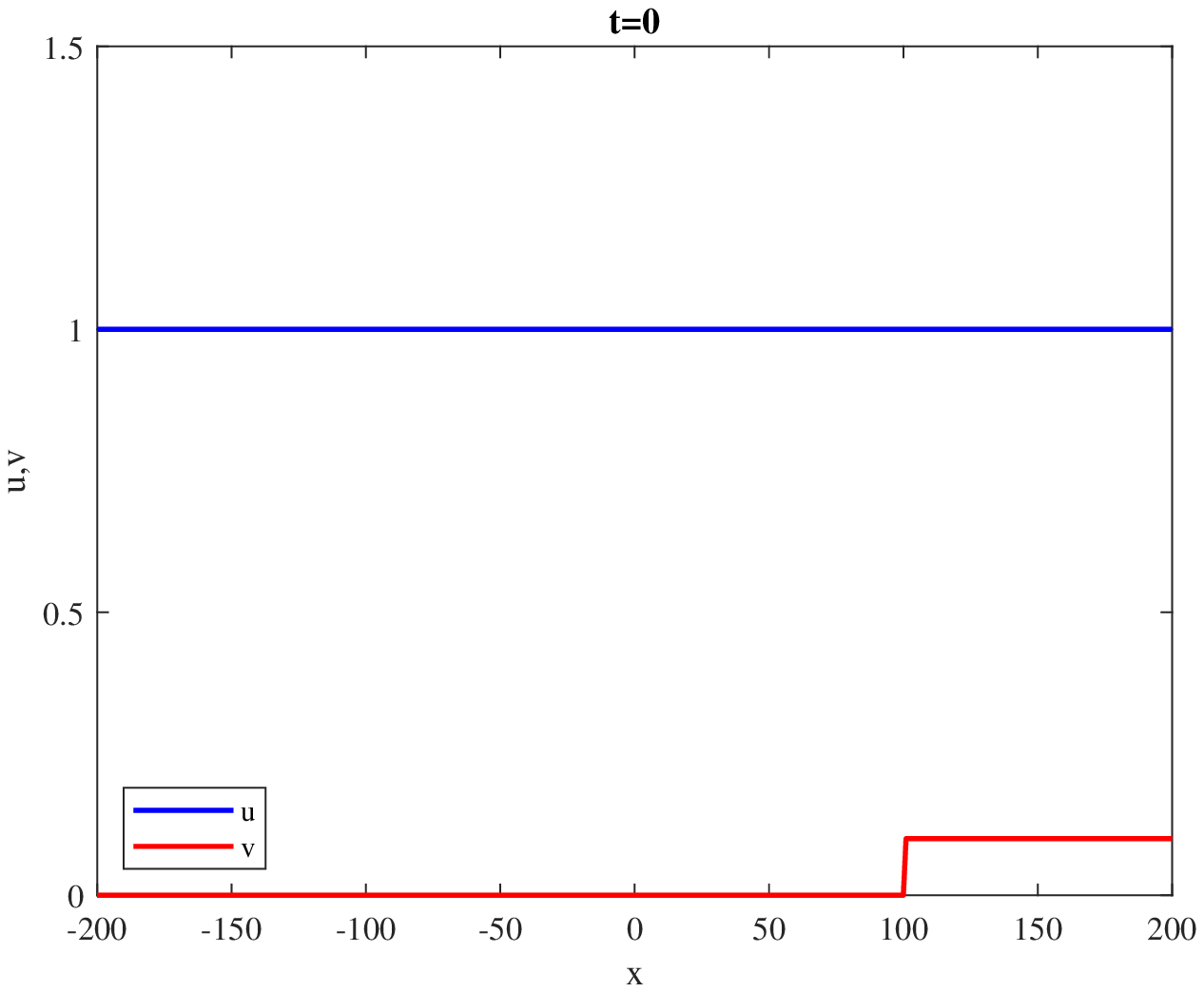}
  \includegraphics[width=2.9in,height=1.7in,clip]{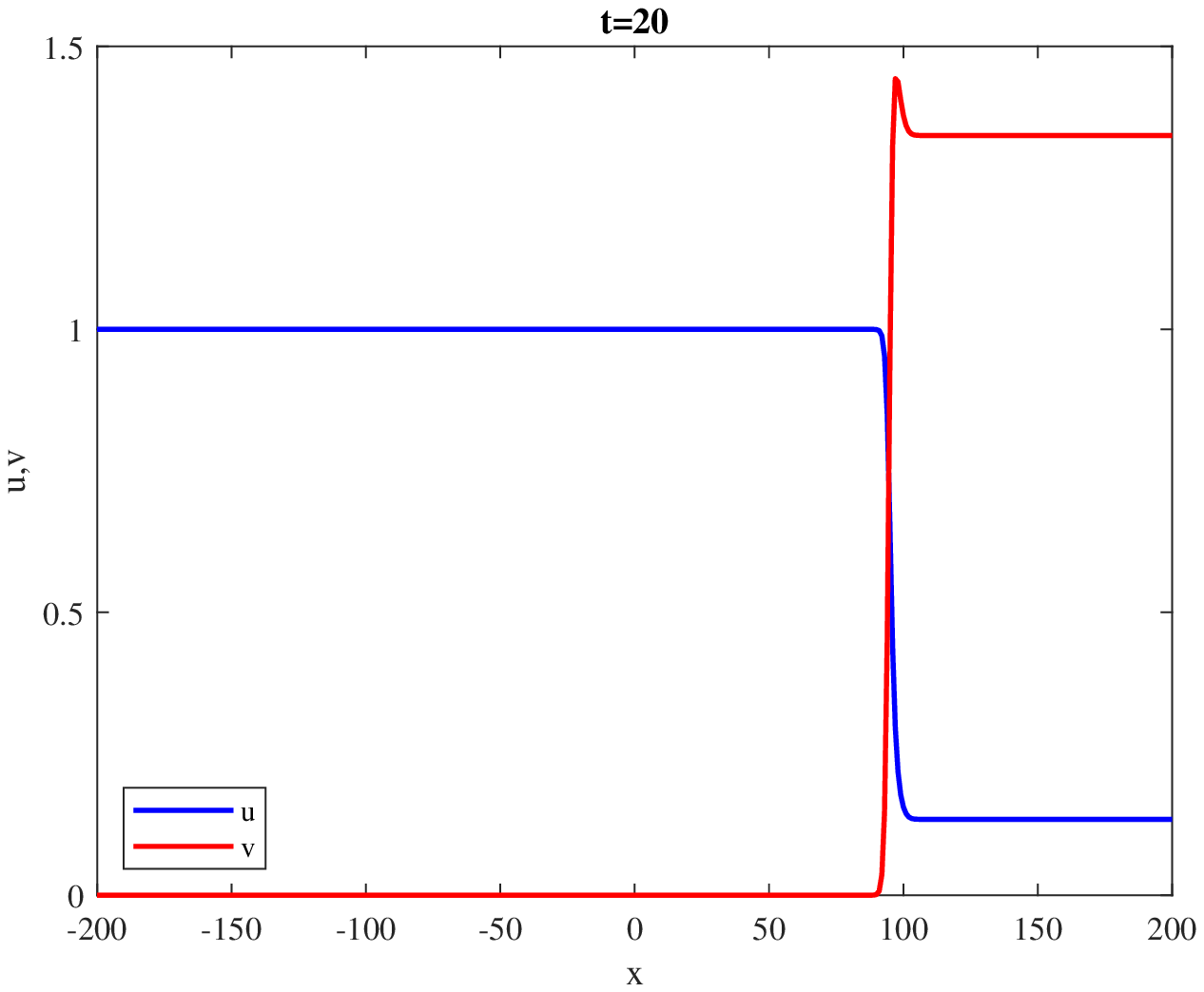}\\
  \includegraphics[width=2.9in,height=1.7in,clip]{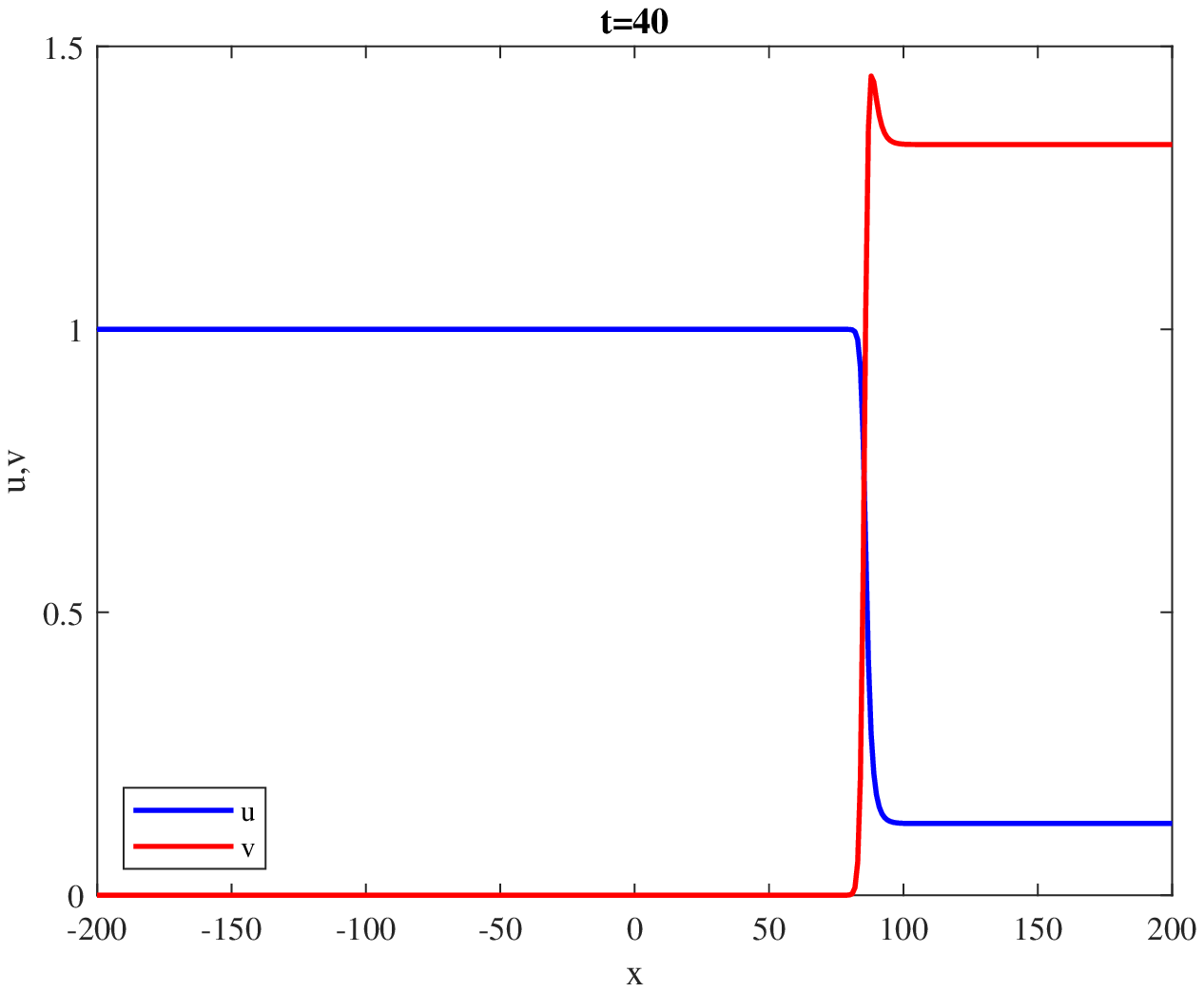}
  \includegraphics[width=2.9in,height=1.7in,clip]{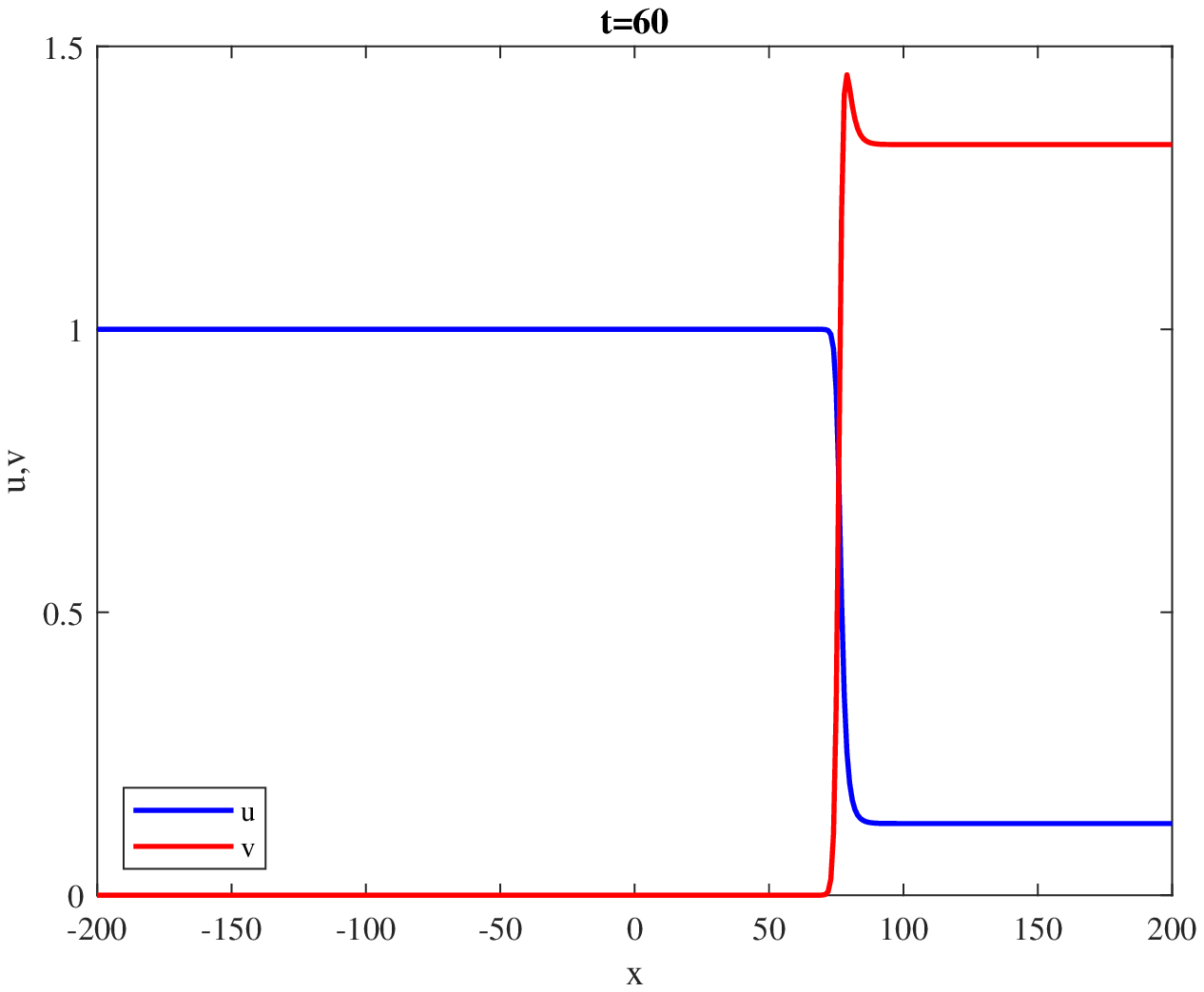}
  \caption{Traveling wave of system \eqref{eq:h2} at different time\\ using  the parameters $d=1$, $s=0.5$, $e_1=2$, $e_2=1.2$ and $a=1.4$.}
  \label{h1}
\end{figure}

\begin{figure}[htbp]
  \centering
  \includegraphics[width=2.9in,height=1.7in,clip]{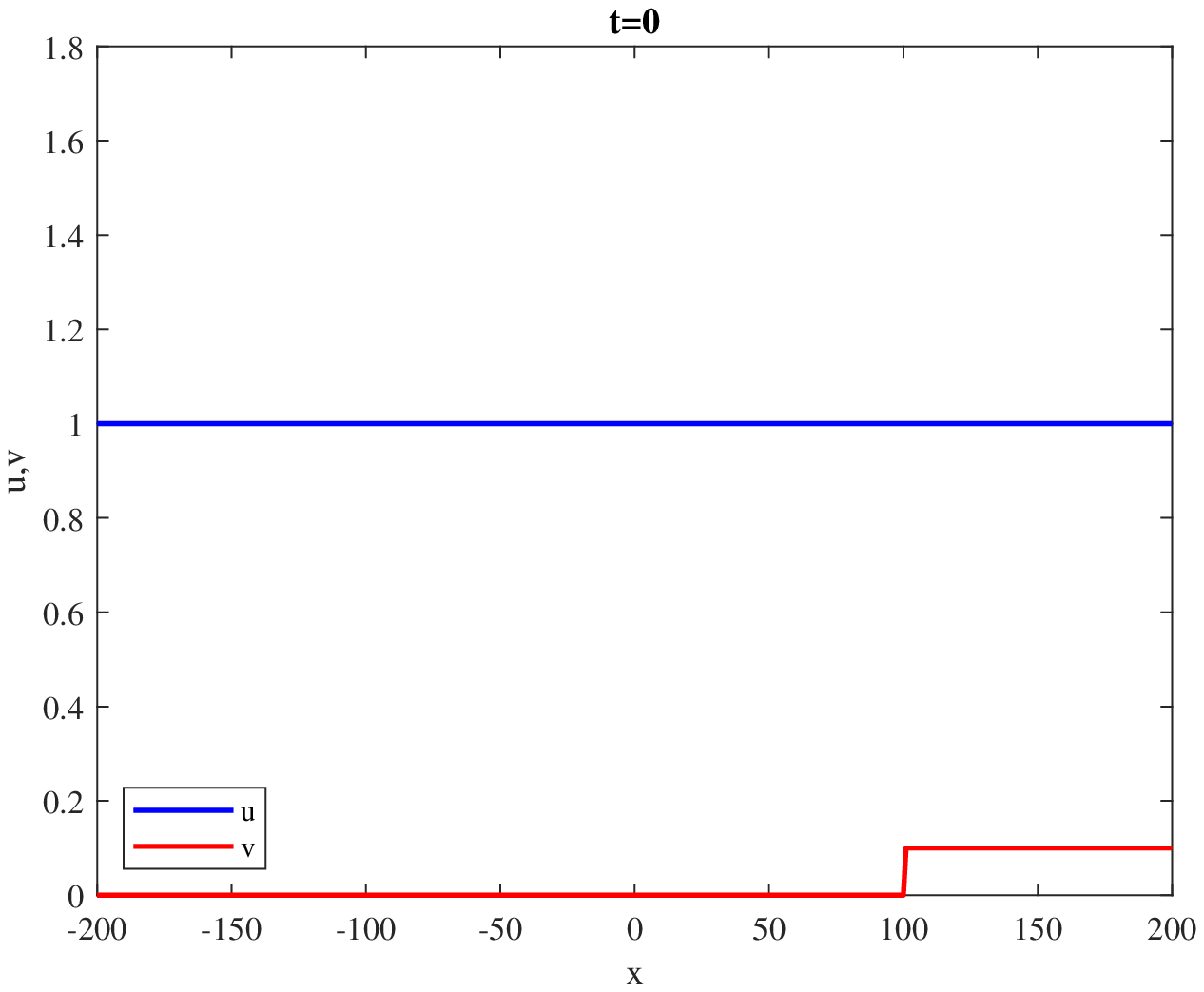}
  \includegraphics[width=2.9in,height=1.7in,clip]{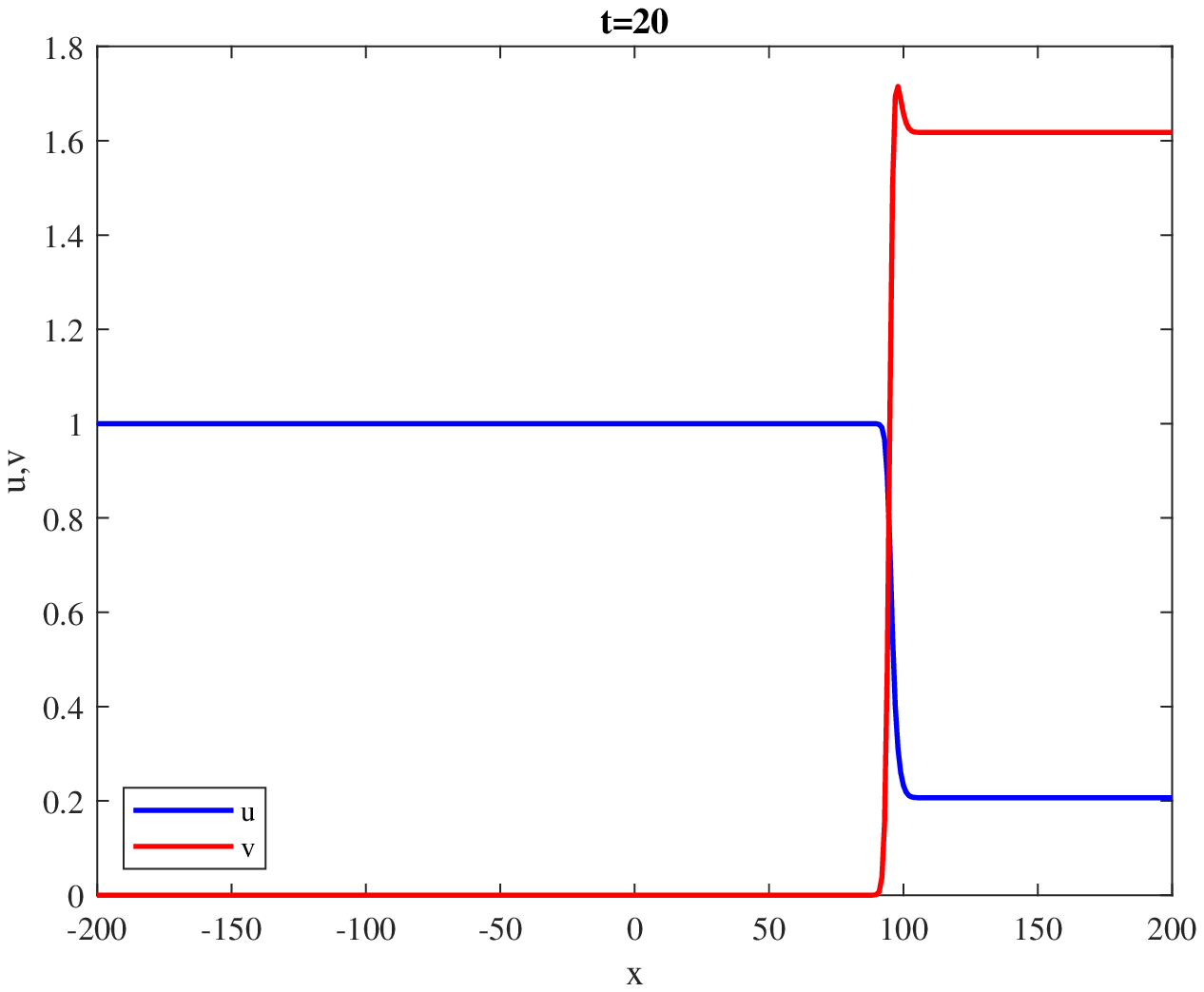}\\
  \includegraphics[width=2.9in,height=1.7in,clip]{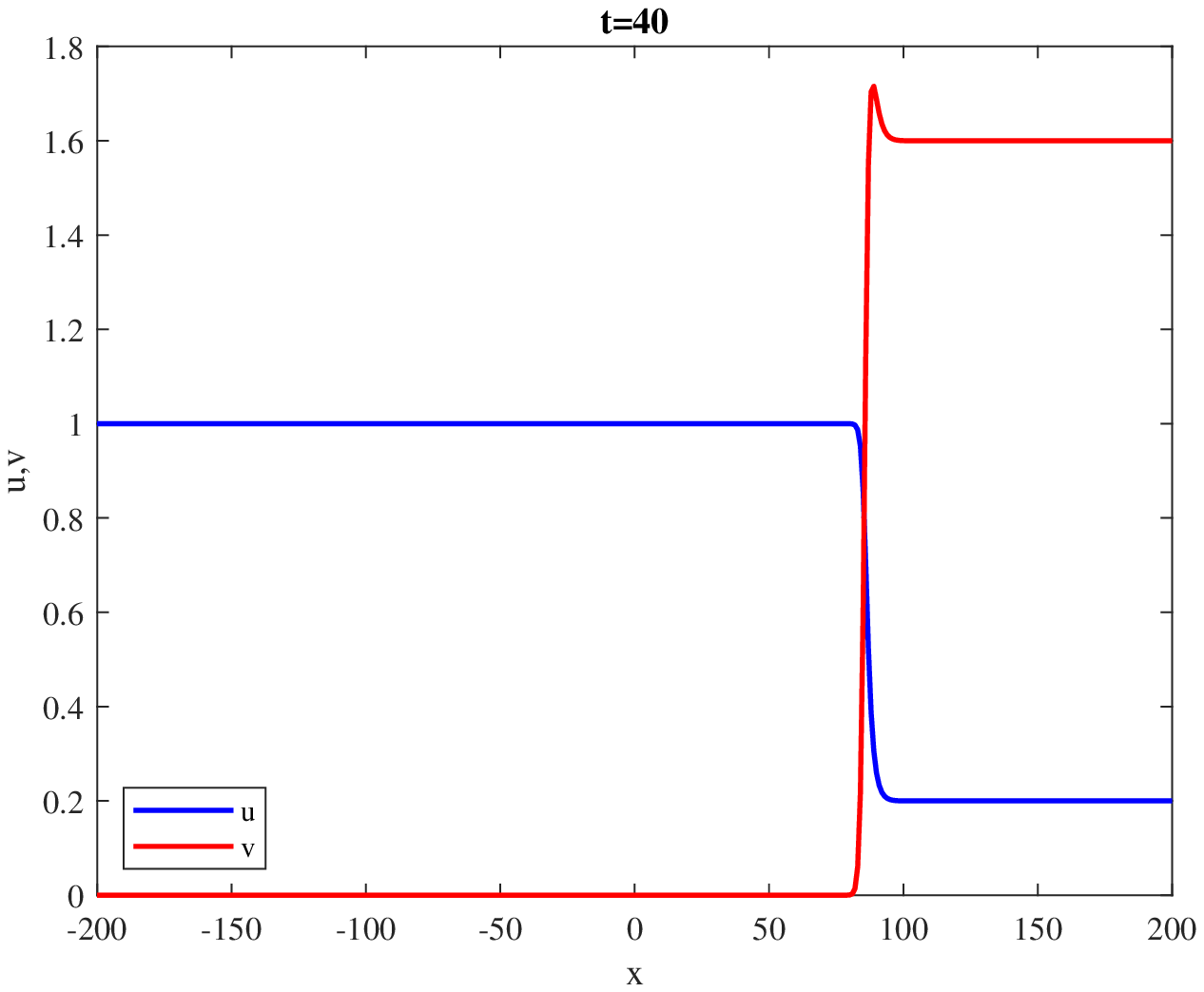}
  \includegraphics[width=2.9in,height=1.7in,clip]{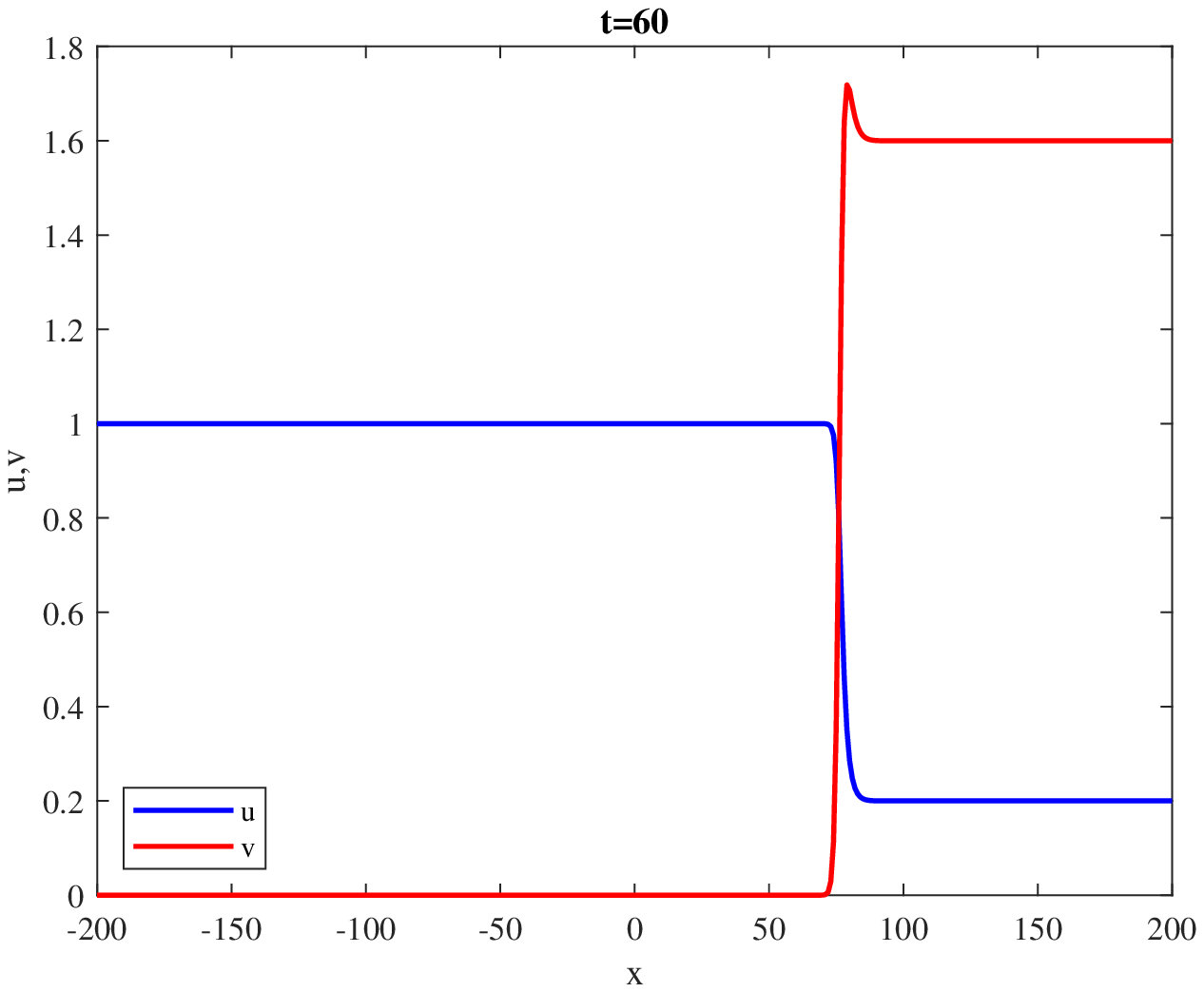}
  \caption{Traveling wave of system \eqref{eq:h2} at different time\\ using  the parameters $d=1$, $s=0.5$, $e_1=1.2$, $e_2=1.4$ and $a=0.7$.}
  \label{h2}
\end{figure}

\begin{figure}[htbp]
  \centering
  \includegraphics[width=2.9in,height=1.7in,clip]{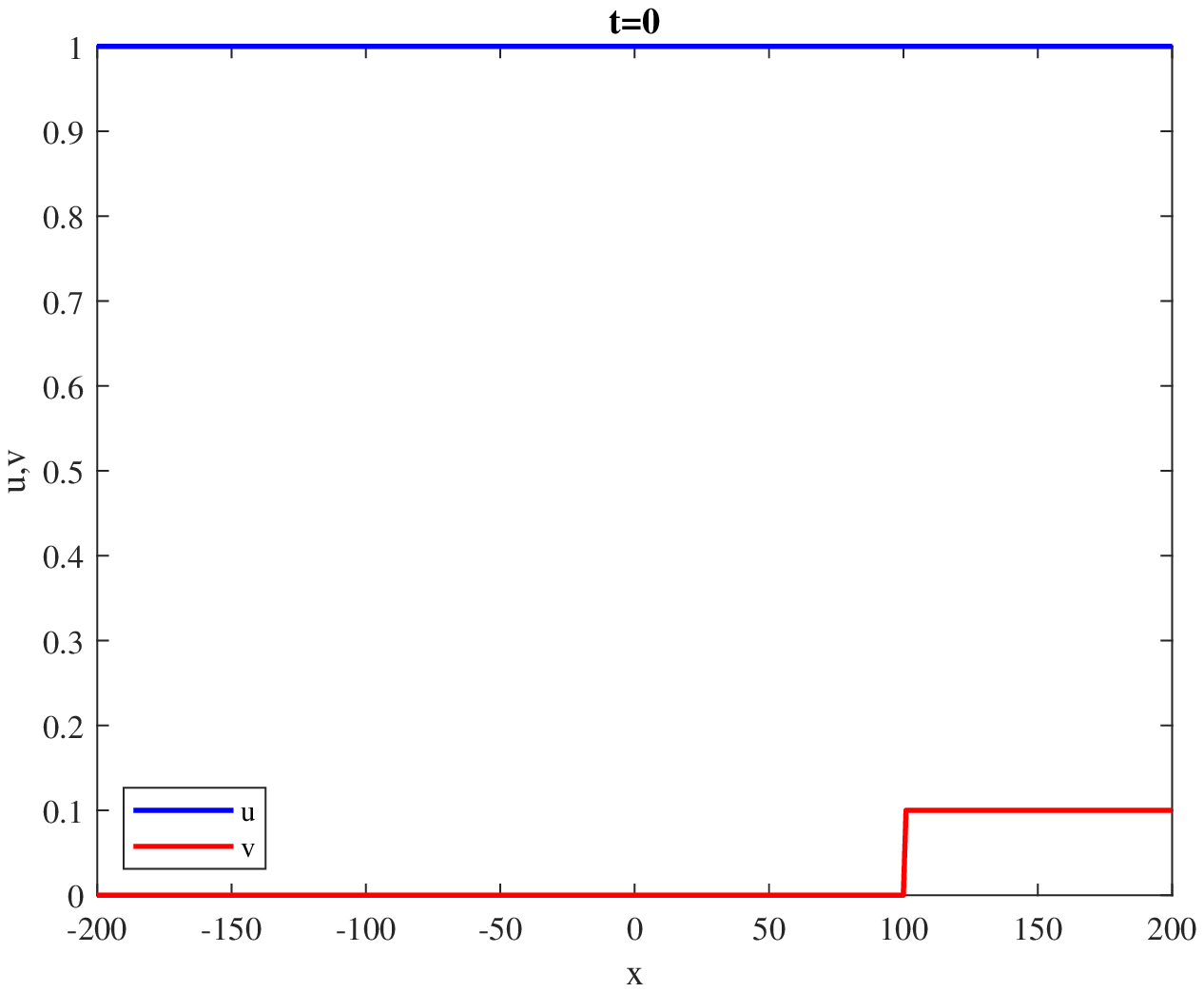}
  \includegraphics[width=2.9in,height=1.7in,clip]{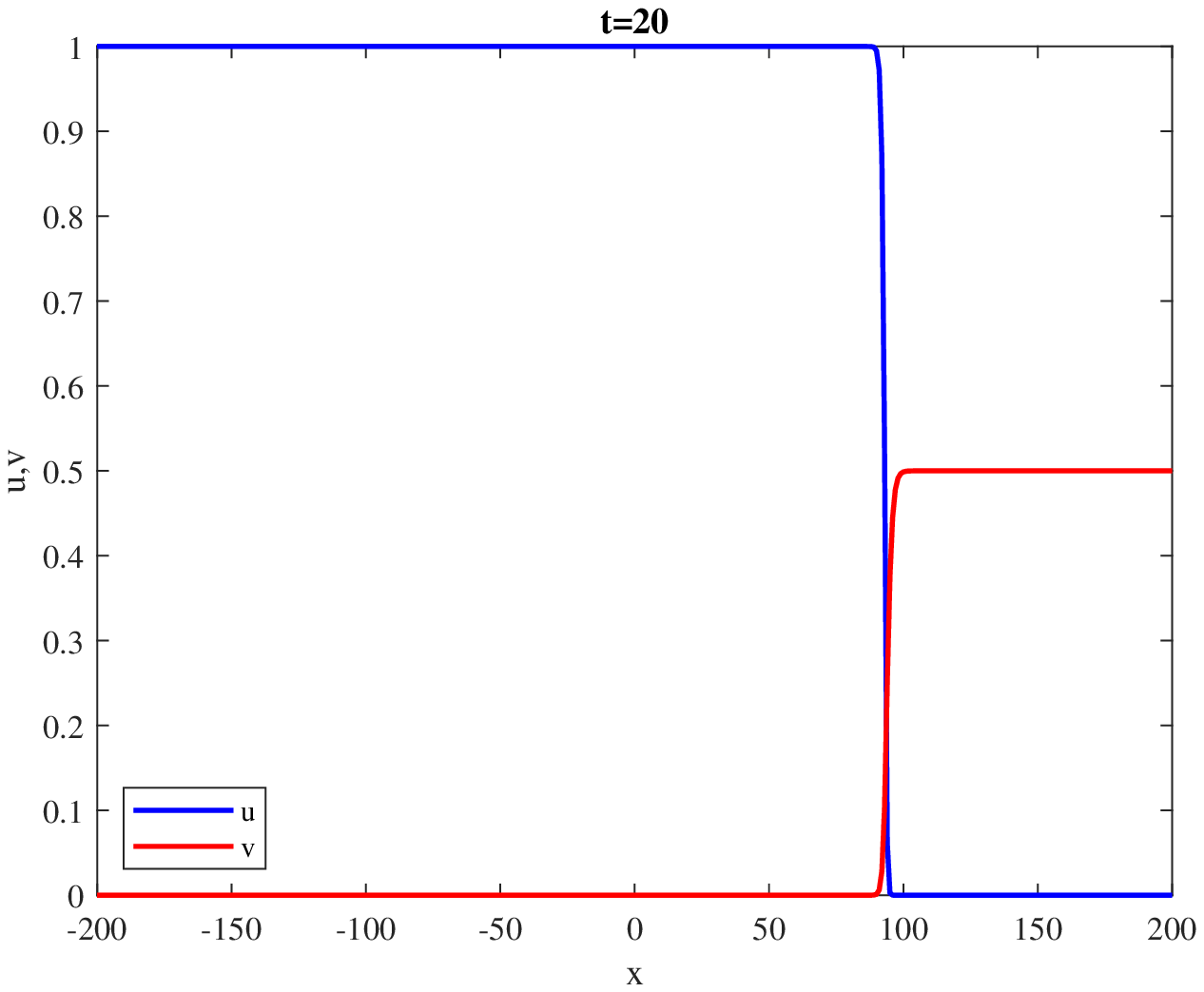}\\
  \includegraphics[width=2.9in,height=1.7in,clip]{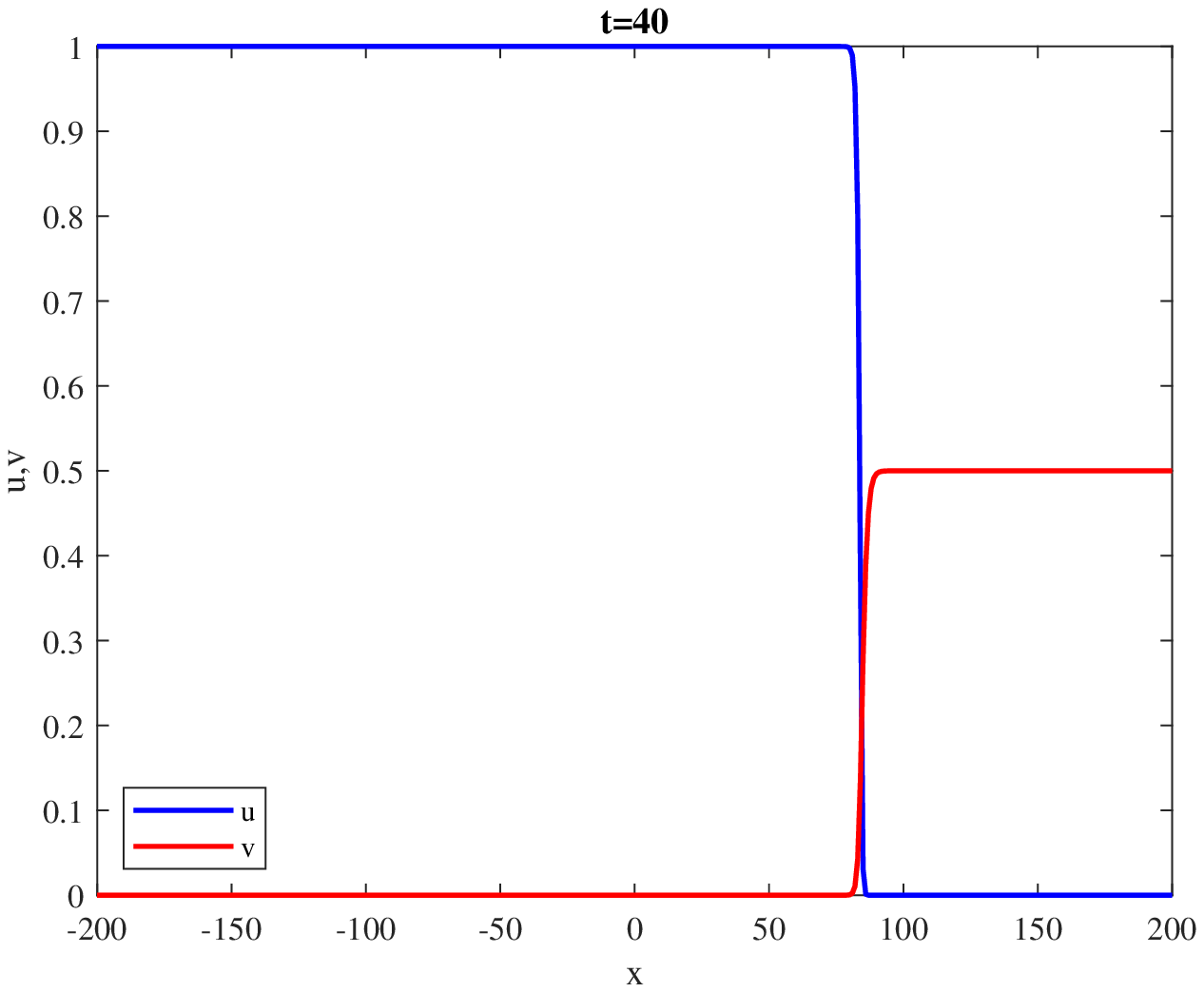}
  \includegraphics[width=2.9in,height=1.7in,clip]{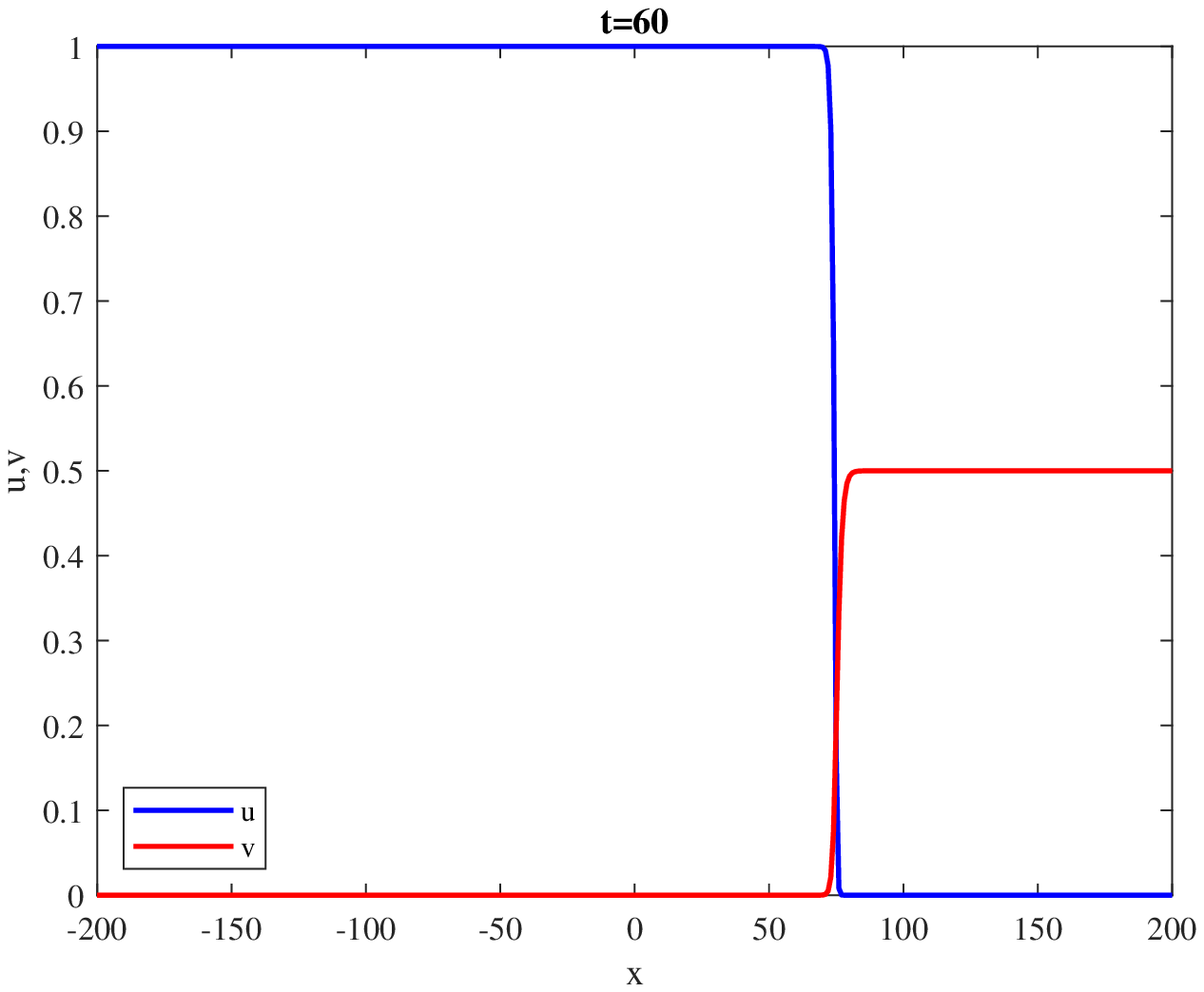}
  \caption{Traveling wave of system \eqref{eq:h2} at different time \\ using the parameters $d=1$, $s=0.5$, $e_1=1.2$, $e_2=0.5$ and $a=15$.}
  \label{h3}
\end{figure}

Below we use MATLAB to further present some numerical simulations. We consider the initial value as follows
\begin{equation}\label{eq:ic}
  u(x,0)=1,  -200\leq x\leq200\ \text{and}\
  v(x,0)=\left\{
  \begin{array}{ll}
    0, &for\ -200<x\leq100,\\[0.2cm]
    0.1,\ &for\ 100<x\leq200.
  \end{array}
  \right.
\end{equation}
First of all, let $d=1$, $s=0.5$, $e_1=2$ and $e_2=1.2$, then $\bar{a}\approx1.4373$. At this time, we set $a=1.4$, then $(u^*,v^*)\approx(0.1266,1.3266)$. On the other hand, let $d=1$, $s=0.5$, $e_1=1.2$ and $e_2=1.4$, then $\bar{a}\approx0.8571$. We further choose $a=0.7$, then $(u^*,v^*)=(0.2,1.6)$. Thus, the above two groups of parameters satisfy different conditions in Theorem \ref{th:gg}. Under the same initial value \eqref{eq:ic}, one can observe traveling wave connecting connecting $(1,0)$ and $(u^*,v^*)$ of system \eqref{eq:h2} from Figure \ref{h1} and Figure \ref{h2}. Moreover, we also select $d=1$, $s=0.5$, $e_1=1.2$, $e_2=0.5$ and $a=15$, which satisfies the condition in Theorem \ref{th:no3}, then system \eqref{eq:h2} has a traveling wave connecting $(1,0)$ and $(0,e_2)$ from Figure \ref{h3}.

\subsection{Modified Leslie-Gower system with Lotka-Volterra type}
\noindent

Let us now consider system
\begin{equation}\label{eq:h1}
\left\{\begin{split}
&\displaystyle\frac{\partial u(x, t)}{\partial t}=\Delta u(x, t)+u(x, t)\left(1-u(x, t)-av(x, t)\right), \\[0.2cm]
&\displaystyle\frac{\partial v(x, t)}{\partial t}=d \Delta v(x, t)+sv(x, t)\left(1-\frac{v(x, t)}{u(x, t)+e_2}\right).
\end{split}
\right.
\end{equation}
where $a$, $d$, $s$ and $e_2$ are positive constants. Then
\begin{align*}
  f(u)=au,  \quad p(u)=\frac{1-u}{a}, \quad q(u)=u+e_2, \quad h(u)=1,\quad \mu=e_2,
\end{align*}
and
\begin{equation*}
  (u^*,v^*)=\left(\frac{1-ae_2}{1+a},\ \frac{1+e_2}{1+a}\right).
\end{equation*}
From Theorem \ref{th:ma2}, one can easily verify that if $a\leq1/(1+e_2)$, then system \eqref{eq:h1} has a traveling wave connecting $(1,0)$ and $(u^*,v^*)$. However, inspired by \cite{D04,wf21}, we obtain a better result by constructing a novel Lyapunov function.

We need the following lemma to proceed.
\begin{lemma}\label{le:sa}
The polynomial $S(a)=a^3+a^2-16a(1+e_2)-16(1+e_2)(2+e_2)$ has a unique positive root $\tilde{a}$ satisfying $\tilde{a}>1/(1+e_2)$.
\end{lemma}
\Proof  Note that the coefficient sign of $S(a)$ changes only once, it follows from Descartes' Rule of Signs that $S(a)$ has a unique positive root $\tilde{a}$. To show that $\tilde{a}>1/(1+e_2)$, we just need to prove $S(1/(1+e_2))<0$. Indeed, we have
\begin{equation*}
  -(1+e_2)^3S\left(\frac{1}{1+e_2}\right)=16(1+e_2)^5+16(1+e_2)^4+16(1+e_2)^3-(1+e_2)-1>0.
\end{equation*}
Hence, we complete the proof.

\qed

With the aid of Lemma \ref{le:sa}, we begin to prove the following theorem.
\begin{theorem}\label{th:3}
For $c\geq c^*$, there is a constant $\tilde{a}>1/(1+e_2)$ satisfying $S(\tilde{a})=0$ such that if $a<\bar{a}=\min\left\{\tilde{a}, 1/e_2\right\}$, then system  \eqref{eq:h1} has a traveling wave connecting $(1,0)$ and $(u^*,v^*)$.
\end{theorem}
\Proof Define
\begin{equation*}
  H(u,v)=\int_{u^*}^u\frac{(\eta+e_2)\eta-(u^*+e_2)u^*}{(\eta+e_2)\eta}d\eta+\varrho\int_{v^*}^v\frac{\eta-v^*}{\eta}d\eta,
\end{equation*}
where $\varrho$ is a positive constant to be determined later. It is clear that
\begin{equation}\label{eq:hd}
\begin{split}
  \partial_u H=&[(u+u^*+e_2)(u-u^*)]/[(u+e_2)u],\quad \partial_v H=\varrho(v-v^*)/v,\\[0.2cm]
  \partial_{uu} H=&\left[u^*(u^*+e_2)(2u+e_2)\right]/[(u+e_2)u]^2,\ \
  \partial_{vv} H=\varrho v^*/v^2.
\end{split}
\end{equation}
Then we define Lyapunov function $\mathcal{L}:\mathbb{R}^4\rightarrow\mathbb{R}$ by
\begin{equation*}
  \mathcal{L}(u,w,v,y)=cH(u,v)-w\partial_u H-dy\partial_v H,
\end{equation*}
which is clearly continuous function with lower bound in $D$.
A simple calculation yields
\begin{equation*}
\begin{split}
  \mathcal{L}'(\chi(z))=&(cw-w')\partial_u H+(cy-dy')\partial_v H-w^2 \partial_{uu} H-y^2 \partial_{vv} H\\[0.2cm]
  =&(u+u^*+e_2)(u-u^*)(1-u-av)/(u+e_2)\\[0.2cm]
  &+\varrho s(u+e_2-v)(v-v^*)/(u+e_2)-w^2 \partial_{uu} H-y^2 \partial_{vv} H\\[0.2cm]
  =&\mathcal{H}(u,v)/(u+e_2)-w^2 \partial_{uu} H-y^2 \partial_{vv} H
\end{split}
\end{equation*}
where
\begin{equation*}
  \mathcal{H}(u,v)=(u+u^*+e_2)(u-u^*)\left(1-u-av\right)+\varrho s(v-v^*)\left(u+e_2-v\right).
\end{equation*}
To show that $\mathcal{L}'(\chi(z))\leq0$ for $u,v>0$, it is sufficient to prove $\mathcal{H}(u,v)\leq0$
owing to $\partial_{uu} H>0$ and $\partial_{vv} H>0$ for $u,v>0$ from \eqref{eq:hd}.
Note that
\begin{align*}
  &1-u-av=-(u-u^*)-a(v-v^*),\\[0.2cm]
  &u+e_2-v=(u-u^*)-(v-v^*),
\end{align*}
by substituting these into $\mathcal{H}(u,v)$, we arrive at
\begin{align*}
  \mathcal{H}(u,v)=&-(u+u^*+e_2)(u-u^*)^2+[\varrho s-a(u+u^*+e_2)](u-u^*)(v-v^*)-\varrho s(v-v^*)^2\\[0.2cm]
  =&-\left[\sqrt{(u+u^*+e_2)}(u-u^*)-\frac{\varrho s-a(u+u^*+e_2)}{2\sqrt{(u+u^*+e_2)}}(v-v^*)\right]^2\\[0.2cm]
  &\quad-\left(\varrho s-\frac{[\varrho s-a(u+u^*+e_2)]^2}{4(u+u^*+e_2)}\right)(v-v^*)^2.
\end{align*}
If $[\varrho s-a(u+u^*+e_2)]^2-4\varrho s(u+u^*+e_2)<0$, that is
\begin{equation}\label{eq:ts}
  R(\varrho s):=(\varrho s)^2-2(a+2)(u+u^*+e_2)\varrho s +a^2(u+u^*+e_2)^2<0,
\end{equation}
then $\mathcal{H}(u,v)\leq0$ for $u,v>0$.
Obviously, $R(x)$ has two roots $0<r_-(u)<r_+(u)$, where
\begin{align*}
  r_-(u)=(u+u^*+e_2)\left(\sqrt{a+1}-1\right)^2,\\[0.2cm]
  r_+(u)=(u+u^*+e_2)\left(\sqrt{a+1}+1\right)^2.
\end{align*}
And \eqref{eq:ts} holds if and only if $r_-(u)<\varrho s<r_+(u)$ for $u\in(0,1)$, which is equivalent to
\begin{equation*}
  \sup_{u\in(0,1)}r_-(u)<\inf_{u\in(0,1)}r_+(u)
\end{equation*}
by the arbitrariness of $u$. More precisely,
\begin{equation*}
  (1+u^*+e_2)\left(\sqrt{a+1}-1\right)^2<(u^*+e_2)\left(\sqrt{a+1}+1\right)^2.
\end{equation*}
Simplifying the above inequality, we have
\begin{equation*}
\left(\sqrt{a+1}-1\right)^2<4(u^*+e_2)\sqrt{a+1}.
\end{equation*}
By using the expression $u^*+e_2=(1+e_2)/(1+a)$, we further get
\begin{equation*}
(a+2)\sqrt{a+1}<2(a+1)+4(1+e_2),
\end{equation*}
which is equivalent to
\begin{equation*}
  S(a)=a^3+a^2-16a(1+e_2)-16(1+e_2)(2+e_2)<0.
\end{equation*}
From Lemma \ref{le:sa}, we complete the proof.

\qed

The last theorem establishes the existence of traveling wave connecting $(1,0)$ and $(0,e_2)$ of system \eqref{eq:h1}.
\begin{theorem}\label{th:no2}
Assume that $c\geq c^*$, if $ae^2_2>1+e_2$, then system \eqref{eq:h1} has a traveling wave connecting $(1,0)$ and $(0,e_2)$ and satisfying $v(z)\in \mathcal{A}\cup \mathcal{B}$ and $w(z)<0$ over $\mathbb{R}$.
\end{theorem}
\Proof Clearly, Assumptions \ref{as} hold and $g(u)=a$. From Theorem \ref{th:ma3}, we complete the proof.

\qed

\begin{figure}[htbp]
  \centering
  \includegraphics[width=2.9in,height=1.7in,clip]{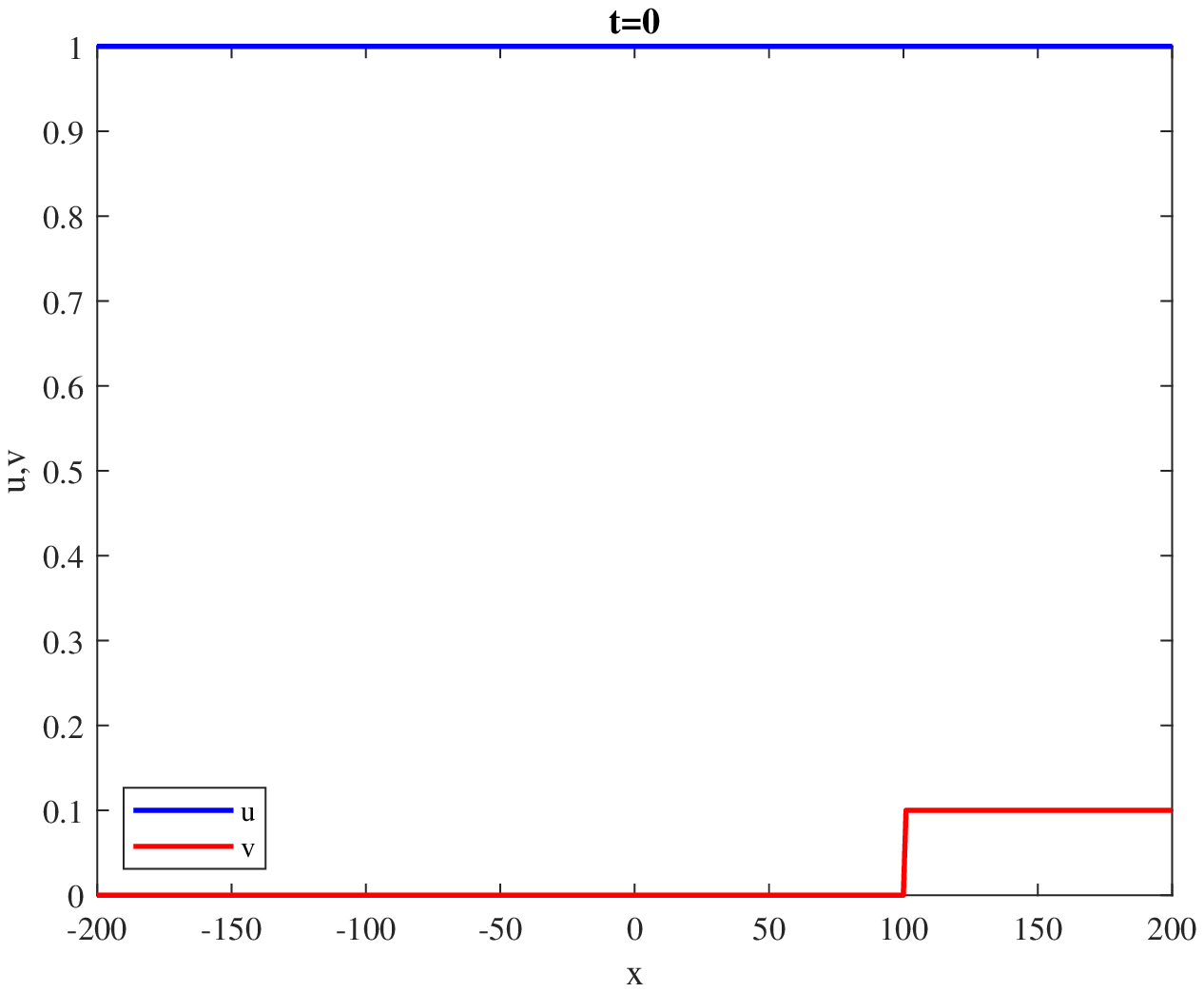}\quad
  \includegraphics[width=2.9in,height=1.7in,clip]{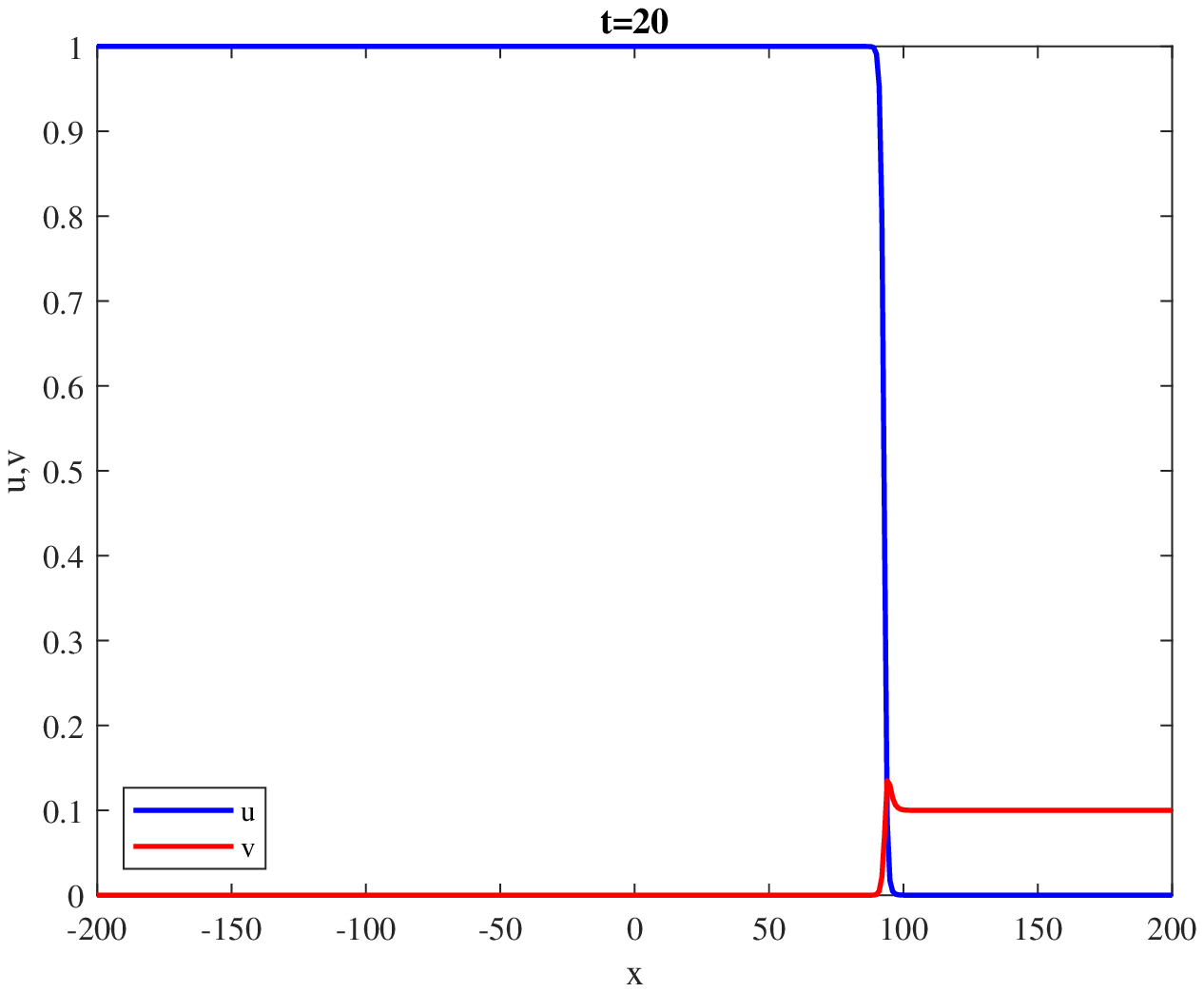}\\
  \includegraphics[width=2.9in,height=1.7in,clip]{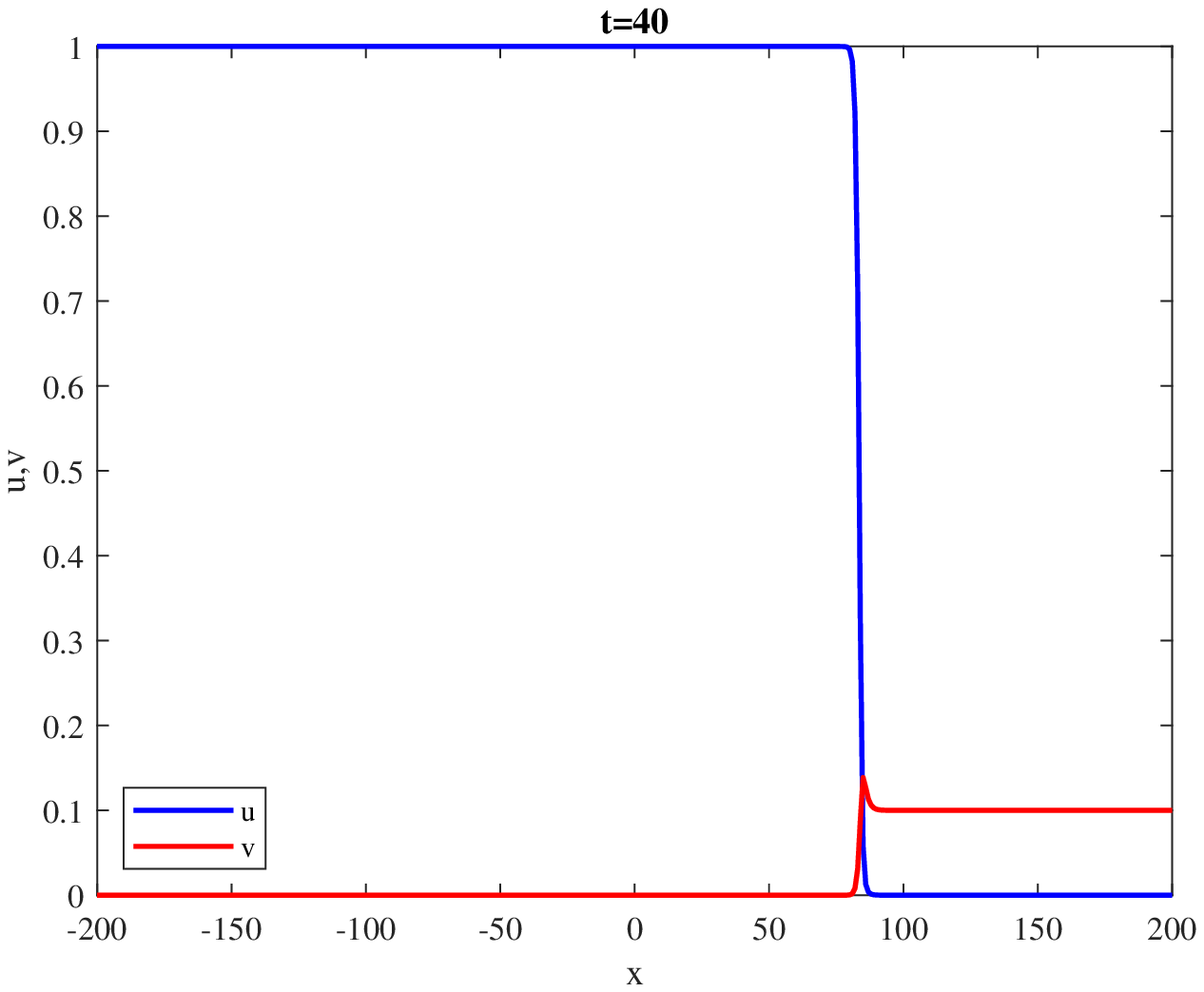}\quad
  \includegraphics[width=2.9in,height=1.7in,clip]{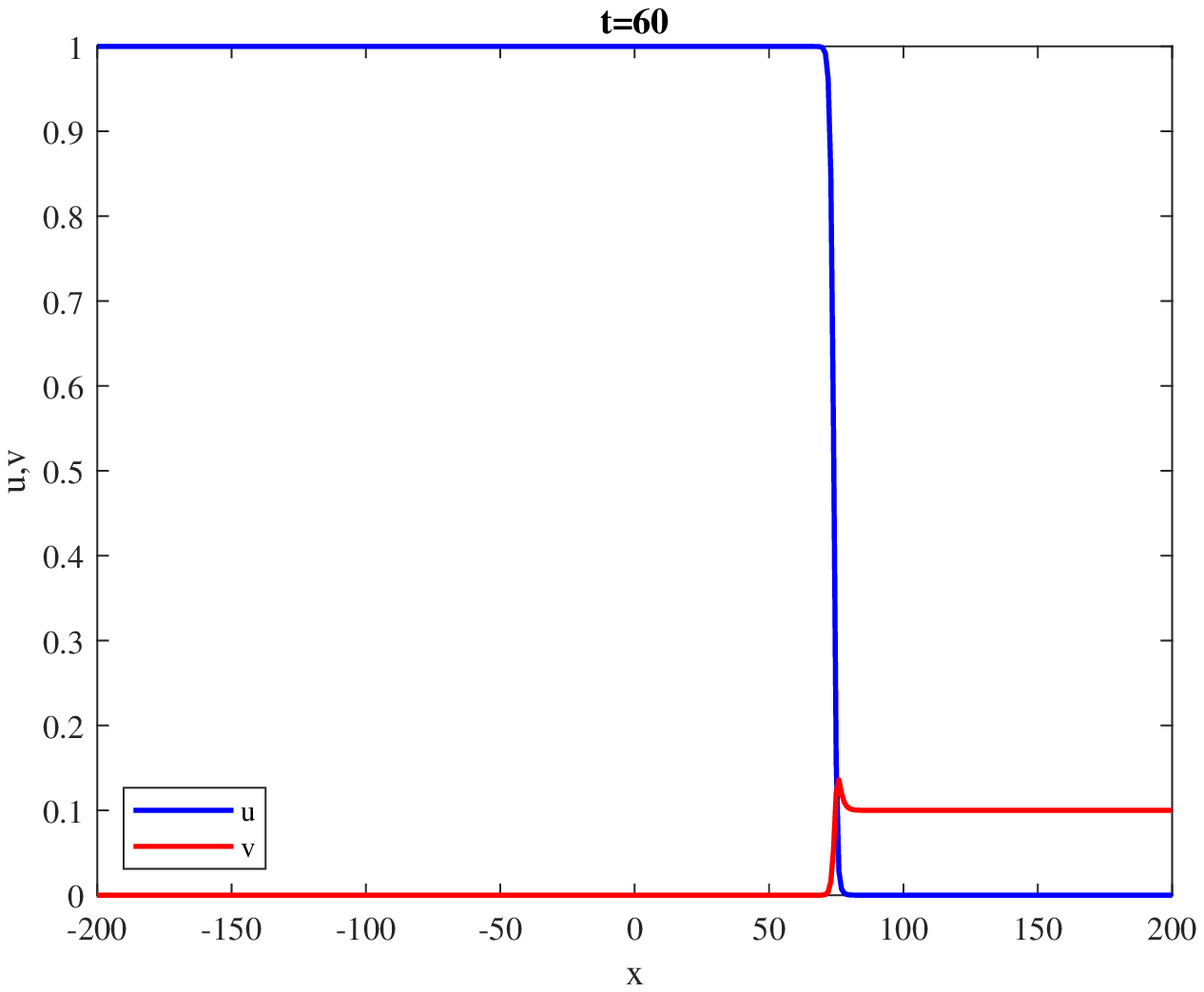}
  \caption{Traveling wave of system \eqref{eq:h2} at different time\\ using  the parameters $d=1$, $s=0.5$, $a=15$ and $e_2=0.1$.}
  \label{c1}
\end{figure}

\begin{figure}[htbp]
  \centering
  \includegraphics[width=2.9in,height=1.7in,clip]{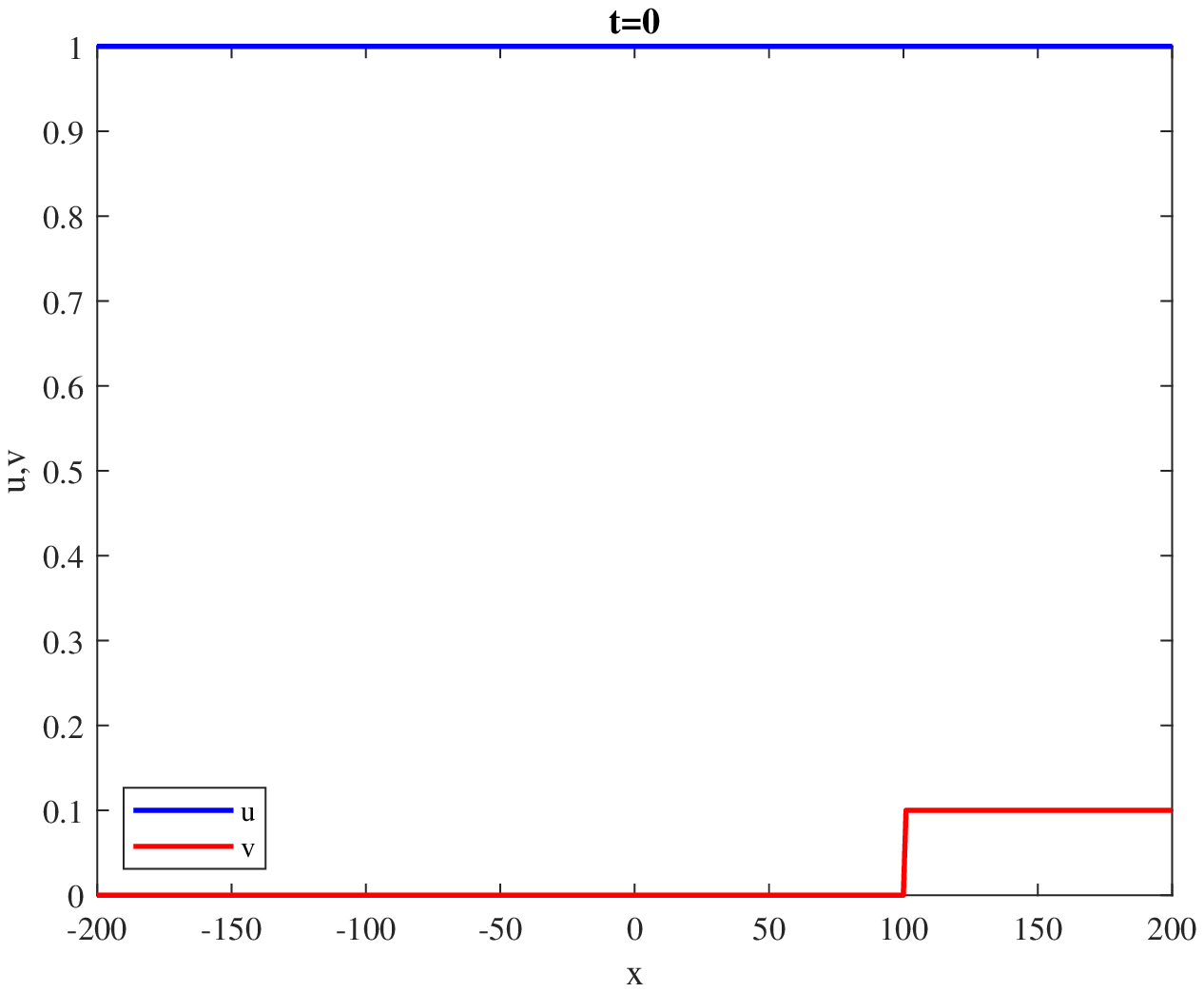}\quad
  \includegraphics[width=2.9in,height=1.7in,clip]{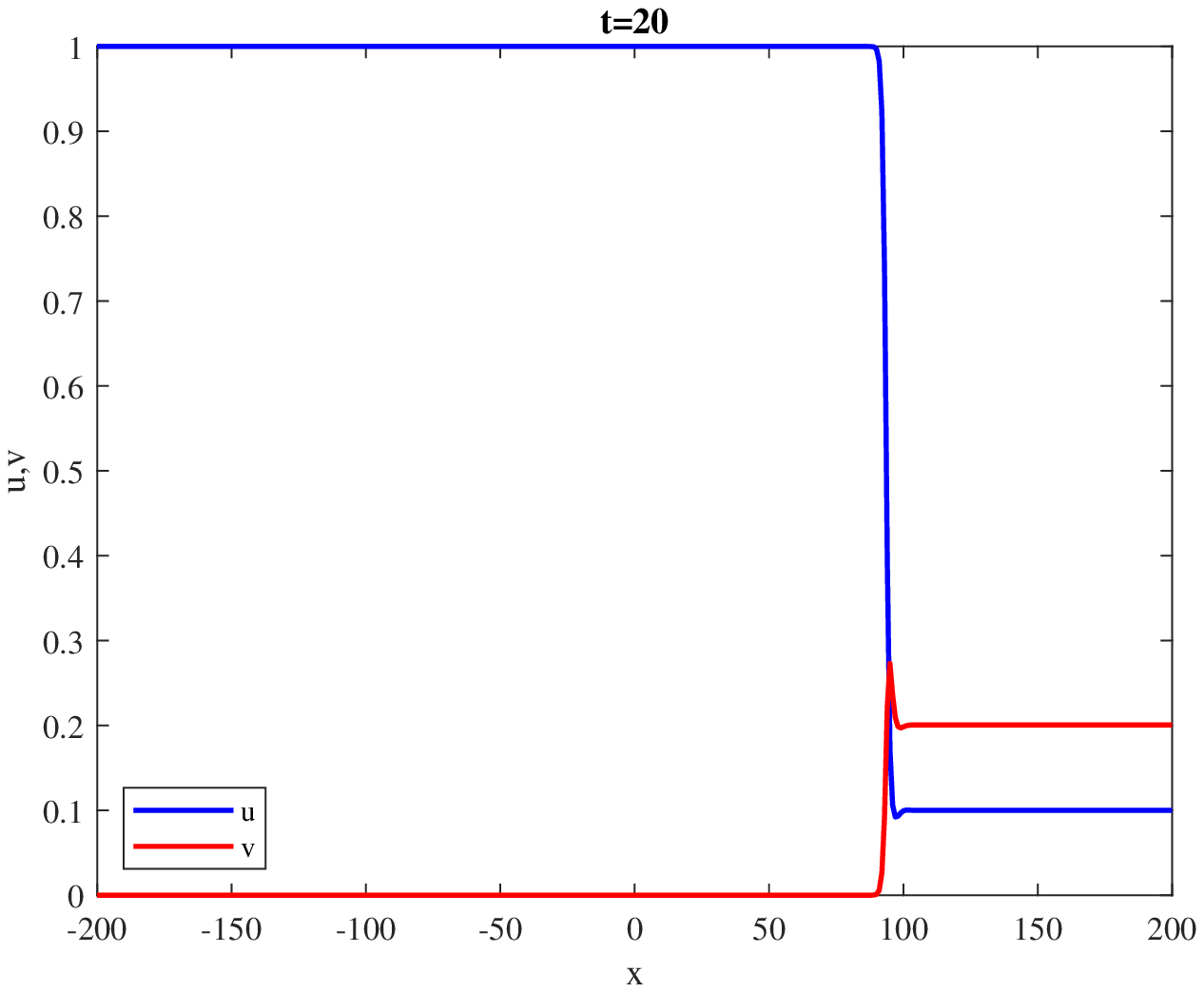}\\
  \includegraphics[width=2.9in,height=1.7in,clip]{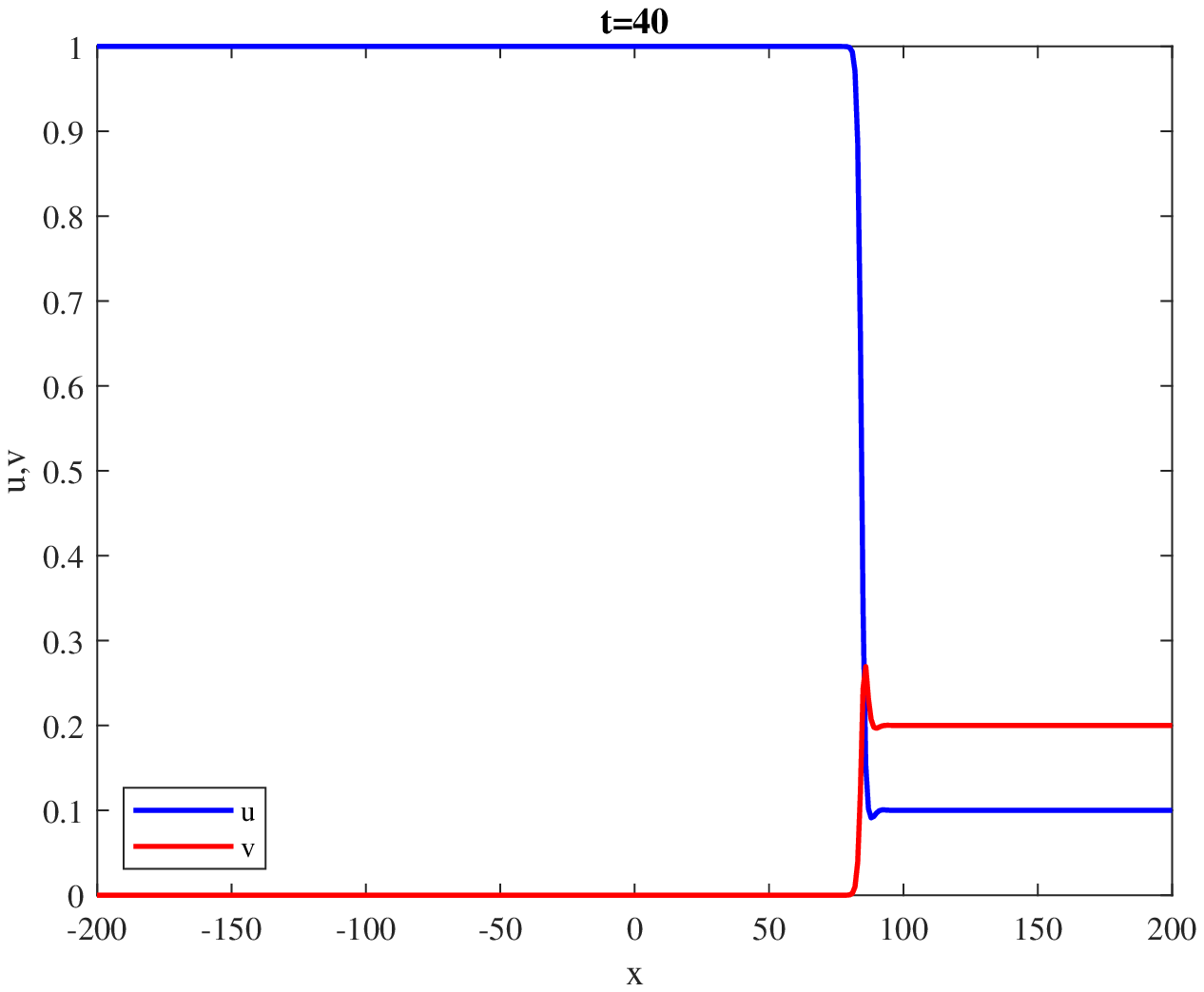}\quad
  \includegraphics[width=2.9in,height=1.7in,clip]{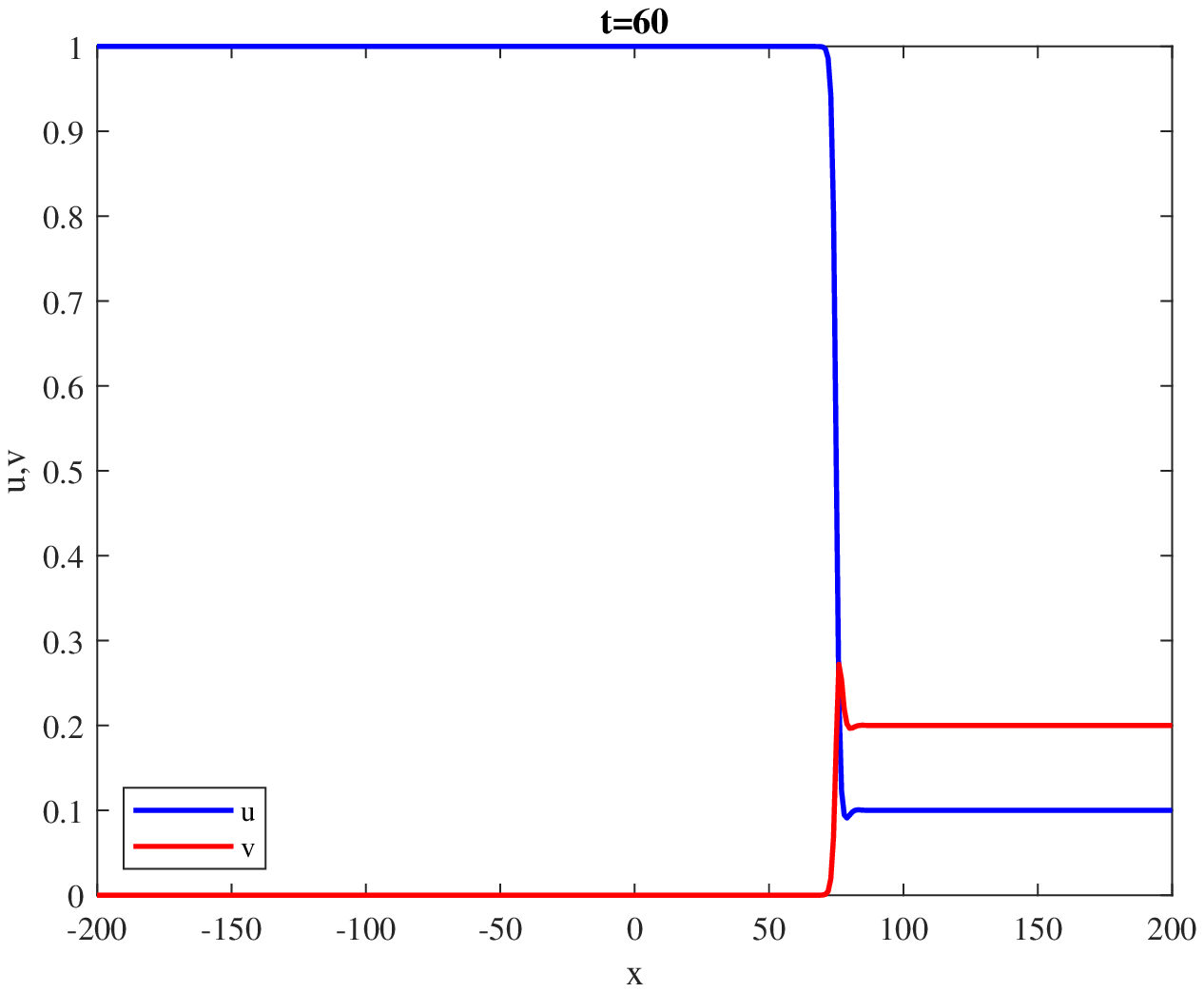}
  \caption{Traveling wave of system \eqref{eq:h1} at different time\\ using  the parameters $d=1$, $s=0.5$, $a =4.5$ and $e_2=0.1$.}
  \label{c2}
\end{figure}

We continue to illustrate our results by numerical simulations under the initial value \eqref{eq:ic}. Let $d=1$, $s=0.5$, $e_2=0.1$ and $a=15$, Figure \ref{c1} depicts that $u$ fades out. Moreover, let $d=1$, $s=0.5$, $e_2=0.1$, then $\bar{a}\approx4.5895$. We further choose $a=4.5$, then $(u^*,v^*)=(0.1,0.2)$. As show in Figure \ref{c2}, system \eqref{eq:h1} has a traveling wave connecting $(1,0)$ and $(u^*,v^*)$.


%

\section*{Acknowledgments}
This work is supported by the National Natural Science Foundation of China (No. 12171039 and 12271044).

\bibliographystyle{elsarticle-num}

\end{document}